\documentstyle[bezier]{article}

\makeatletter

\newcounter{oekaki}                
\def\theoekaki{\@arabic\c@oekaki}  
\def\fps@oekaki{hbt}
\def\ftype@oekaki{8}
\def\ext@oekaki{loe}
\def\fnum@oekaki{Fig. \theoekaki}
\def\oekaki{\@float{oekaki}}
\let\endoekaki\end@float
\long\def\@makecaption#1{          
  \setbox\@tempboxa\hbox{[ #1 ]}
  \ifdim \wd\@tempboxa >\hsize
   #1 \par
  \else
   \hbox to\hsize{\hfil\box\@tempboxa\hfil}
  \fi}

\makeatother
\newtheorem{prop}{Proposition}

\newtheorem{lem}{Lemma}
\newtheorem{thm}{Theorem}

\newcommand{\np}[0]{\mbox{\bf :}}               
\newcommand{\fb}[0]{\raisebox{.6ex}{\framebox[0.5em]{}}}
\newcommand{\itg}[0]{\mbox{\bf Z}}              
\newcommand{\cpx}[0]{\mbox{\bf C}}              
\newcommand{\og}[4]{{#1}^{#4}_{#3 #2}}          
\newcommand{\iv}[3]{\bar{\phi}({#1})^{#2}_{#3}{}} 
\newcommand{\ov}[3]{\phi({#1})^{#2}_{#3}{}}     
\newcommand{\path}[3]{{\cal P}({#1})^{#2}_{#3}} 

\newcommand{\tr}[0]{{\rm Tr}}              

\newcommand{\PR}[1]{\check{R}({#1})}            
\newcommand{\PW}[1]{\check{W}({#1})}             
\newcommand{\pw}[0]{\check{W}}                  

\begin{document}


\title
{\LARGE
{\bf
Ruijsenaars' commuting difference operators\
as commuting transfer matrices
\footnote
{Short title : Ruijsenaars' operators as commuting transfer matrices.
AMS subject classification: 81R50, 17B37, 05E05.}
}}

\author{Koji HASEGAWA
\thanks
{Partly supported by
The Kawai Foundation for the Promotion of Mathematical Science,
Japan Society for the Promotion of Science,
and Grant-in-Aid for Scientific Research
on Priority Areas, the Ministry of Education, Science, Sports
and Culture, Japan.}
\\Mathematical Institute, Tohoku University, Sendai JAPAN}

\date{}
\maketitle

\begin{abstract}
\noindent
For Belavin's elliptic quantum R-matrix, we construct an L-operator
as a set of difference operators acting on functions on
the type A weight space. According to the fundamental relation
$RLL=LLR$, the trace of the L-operator gives a
commuting difference operators.
We show that for the above mentioned L-operator this approach gives
Macdonald type operators with elliptic theta function coefficient,
actually equivalent to Ruijsenaars' operators.
The relationship between the difference L-operator and Krichever's
Lax matrix is given, and an explicit formula for elliptic commuting
differential operators is derived.
We also study the invariant subspace for the system
which is spanned by symmetric theta functions
on the weight space.

\end{abstract}

\section{Introduction}

In \cite{Mac88}, \cite{Mac90}, I.~G.~Macdonald defined a commuting
system of difference operators for each root system
and a new family of orthogonal polynomials
containing two rational parameters $(q,t)$
(in case all the roots have the equal length)
as their similteneous eigenfunctions.
Up to now,
there are at least two theories which provide understanding
of the system with general value of  $(q,t)$.
One is the work by Etingof and Kirillov,
who obtained these operators as the image of
central elements of the quantum enveloping algebra $U_q({sl}_n)$
acting on ``vector valued characters'' \cite{EK1}.
The other is the work by Cherednik \cite{C92},
who used double affine Hecke algebra,
its representation via $q$-difference operators,
and the center of the algebra.
(See also Section 7.)

Here we wish to suggest yet another approach for the system.

Needless to say, the Yang-Baxter equation is one of the
important backgrounds of the above two works.
Originally,
in Baxter's study of two-dimensional lattice statistical models,
the Yang-Baxter equation arose as the condition
to provide sufficiently many commuting operators.
This is done by taking the trace of the so-called L-operators,
the operators which satisfy the
``$RLL=LLR$  relation" (\ref{eq:RLL=LLR}).
In the lattice model situation, the R-matrix reads as the local
Boltzmann weight of the model.
Given an R-matrix satisfying the Yang-Baxter equation,
an L-operator naturally arises as the row-to-row transfer
matrix and then its trace gives
the commuting transfer matrix \cite{Bax71}\cite{xyz}.

Apart from the lattice models,
recall the Lax matrix method in the completely integrable systems.
Once the equation in problem is formulated in a Lax form,
natural candidates of commuting integrals of the motion
are the characteristic polynomials of the Lax matrix,
namely the traces of the powers of the
matrix.
Their commutativity does not hold in general but one can make use
of the r-matrix structure, a differential form of the
$RLL=LLR$ relation, to ensure it
(see \cite{Skldr} and the references therein).

Therefore the following question seems quite natural to ask:
``what kind of operator arise if we start with the L-operator
realized as difference operators for appropriate functions and take
the trace?"
This is our approach and we will show that
this idea actually works quite well at least for one interesting case.

The case we consider in this paper is
for the elliptic R-matrix of Belavin
\cite{Belavin}
and the main goal is the reproduction of Ruijesenaars' elliptic system.
In the trigonometric limit,
up to a certain simple ``gauge transformation'' \cite{Rgauge} ,
this R-matrix degenerates to the image of
the universal R-matrix for the quantum affine enveloping algebra
$U_q(A_{n-1}^{(1)})$ (\cite{J},\cite{Dr})
in the vector representation.
The $n=2$ case of Belavin's R-matrix is
nothing but Baxter's eight vertex model, for which the
novel quadratic algebra is defined by Sklyanin \cite{Sklyanin}.
We hope that the present paper will give some insight into the
structure related to the Sklyanin algebra and Belavin's elliptic R-matrix.

\medskip
The plan of this paper is as follows.
For Belavin's $n$-state R-matrix, we have constructed an L-operator whose
matrix elements are a certain difference operators
in the previous paper \cite {has}(Section 2, Theorem 1).
This is an $sl_n$- generalization of the L-operator corresponding to the
Sklyanin's difference operators in \cite{Sklyanin}.
For this L-operator, we apply the fusion procedure (Section 3)
and compute the traces of the fused L-operator (Section 4).
The $RLL=LLR$ relation ensures the commutativity of the resulting
operators.
It turns out that the traces exactly give us Macdonald-type
difference operators $M_d$ ($d=1,\cdots,n$), whose
coefficients are given in terms of
Jacobi's elliptic theta function (Theorem 2, main result).
Fixing the elliptic modulus parameter $\tau$,
Belavin's R-matrix has one parameter $\hbar$, while the factorized
L-operator admits another parameter $c$ as well as the spectral
parameter $u$.
Consequently, the operator $M_d$ depends on these three parameters
$\hbar,c,u$
and it turns out that the spectral parameter $u$ appears
only in the overall factor.
The remaining two parameters just play the role of Macdonald's parameters
$q,t$ and the explicit correspondance is given by
$q=\exp \pi \sqrt{-1}\hbar$, $t= q^{-\frac{c}{n}}$.

The computation of the trace of the fused L-operators hinges upon an
interesting theta function determinant formula (Lemma 1).
It contains the parameter $\hbar$.
In the limit $\hbar\rightarrow 0$,
this formula degenerates to the well-known Cauchy type determinant,
which probably first appeared as formula (12) of \cite{Frob}
in the literature, also known as the genus one case of
Fay's trisecant formula \cite{Fay}.

In Section 5, we give explicit relationship between our elliptic
commuting operators and other approaches for the system.
First of all, consider the trigonometric limit.
Then it is easy to see that
our operators tend to Macdonald's operators.
In this trigonometric case another presentation of the system
is known, that is, the trigonometric Ruijsenaars' operators.
The equivalence of Macdonald's operators and Ruijsenaars' ones
is given via the conjugation by a certain function multiplication
(due to T.Koornwinder; see \cite{D}).
In Subsection 5.1 (Proposition 2) we extend this equivalence
to the present elliptic case, namely
between our commuting difference operators and
Ruijsenaars' relativistic elliptic Calogero-Moser system
\cite{Ruij}.
Secondly, we look at the differential limit case, i.e.
the elliptic Calogero-Moser system.
This system admits the Lax formalism 
\cite{Krichever}.
Now the question is whether and how Krichever's Lax matrix arise from
the factorized L-operator;
this is answered in Subsection 5.2 (Proposition 3).
Subsection 5.3 is devoted to the generating function for the operators
$\{M_d\}$.
Since $M_d$ is defined to be the trace of the factorized
L-operator in the ``degree $d$ exterior'' fused representation,
we expect that the generating function is
just the characteristic polynomial of the L-operator.
This is in fact the case: Theorem 3.
The formula can be considered as an elliptic extension of the
generating function in Macdonald's case (\cite{Macbook},VI(3.2)).
An analogue of Jiro Sekiguchi's generating
function \cite{JiroSekiguchi} is obtained as well.
Further, it turns out that our formula is useful to derive
the commuting differential operators in an explicit way
 : we will give an elliptic generalization of
A.Debiard's \cite{De} formula.

As in the theory of Macdonald polynomials,
one may consider the family of symmetric functions
defined as the joint eigenfunction for $\{M_d\}_{d=1,\cdots,n}$.
This diagonalization problem is still under investigation.
Here we would like to establish some structure theorem on a certain
invariant subspace for the L-operator.
Setting the parameter $c$ to be a nonnegative integer $l$,
it is the space $Th_l^{S(n)}$ spanned by the level $l$
$A_{n-1}^{(1)}$ affine Lie algebra characters,
or the symmetric theta functions on the weight space of type A.
This space is a higher rank analogue of Sklyanin's
finite dimensional function space.
We will state that this space can be identified with the
symmetrically fused representation
as the module of the L-operator algebra (Theorem 5).

Brief discussion is given in Section 7.
Appendix is devoted to proving the theta function identity Lemma 1.

\label{sec:intro}


\section{Review of the factorized L operator}

For $n>1$ let
$
V
= \oplus_{k\in \itg / n\itg} {\cpx}e^k
$
(
$e^k=e^{k+n}$)
 and let
$g,h \in {\rm GL}(V)$
 be
$
g e^k := e^k {\rm exp}\frac{2\pi ik}{n}
,
h e^k := e^{k+1}.
$
We have
$
gh = hg{\rm exp} \frac{2\pi i}{n}.
$
Let ${\hbar},\tau \in {\cpx},  {\rm Im}\tau >0$.
We assume
$\hbar \notin {\bf Z} + {\bf Z}\tau$ for convenience.

Belavin's R-matrix $R(u)=R_{\hbar}(u)$
is characterized as the unique solution of the following five conditions.
\begin{eqnarray*}
&&
\hspace{-6em}
\bullet         
\mbox{
        $R_\hbar(u)$ is a holomorphic
        ${\rm End}(V \otimes V)$
        -valued function in $u$,
        }
\label{characterizeR0}
\\
&&
\hspace{-6em}
\bullet         
\mbox{
        $R_\hbar(u)= (x \otimes x)R_\hbar(u)(x \otimes x)^{-1}$
        for $x = g, h$,
        }
\label{characterizeR1}
\\
&&
\hspace{-6em}
\bullet         
\mbox{
        $R_\hbar(u+1)=
        (g\otimes 1)^{-1} R_\hbar(u) (g \otimes 1) \times (-1)$,
        }
\label{characterizeR2}
\\
&&
\hspace{-6em}
\bullet         
\mbox{
        $R_\hbar(u+\tau)=
        (h\otimes 1) R_\hbar(u) (h \otimes 1)^{-1} \times
                (-{\rm exp}2\pi i(u+\frac{\hbar}{n}+\frac{\tau}{2}))^{-1}$,
        }
\label{characterizeR3}
\\
&&
\hspace{-6em}
\bullet         
\mbox{
        $R_\hbar(0)=P:x\otimes y \mapsto y\otimes x$.
        }
\label{characterizeR4}
\end{eqnarray*}
\noindent
One verifies that
 1)\,
there is a unique solution to the above conditions and
 2)\,
the solution satisfies the Yang-Baxter equation.
An explicit formula for $R(u)$ is also available.
Put
\begin{equation}
\theta_{m,l}(u,\tau) :=
        \sum_{\mu \in m+l{\bf Z}}^{}
        {\exp 2\pi i(\mu u + \frac{\mu^2}{2l}\tau)}
\label{eq:theta_jl}
\end{equation}
and
$
\theta^{(j)}(u):=
\theta_{\frac{1}{2}-\frac{j}{n},1}(u+\frac{1}{2}, n\tau).
$
The zeroes of $\theta^{(j)}(u)$ are given by
$\itg + (j+n\itg)\tau.$
Then we have \cite{RT}
$$
R(u) e^i \otimes e^j
=
\sum_{i',j'=1}^n e^{i'} \otimes e^{j'} R(u) ^{ij}_{i'j'},
$$
\begin{equation}
R(u) ^{ij}_{i'j'} =
\delta_{i+j, i'+j' {\rm mod} n}
\frac
        {\theta^{(i'-j')} (u+{\hbar})}
        {\theta^{(i'-i)} ({\hbar})\theta^{(i-j')}(u)}
\frac
        {\prod_{k=0}^{n-1}{\theta^{(k)}(u)}}
        {\prod_{k=1}^{n-1}{\theta^{(k)}(0)}}
{}.
\label{eq:RTformula}
\end{equation}
See figure 1 for our convention to visualize the matrix $R$.
\begin{oekaki}\refstepcounter{oekaki}\addtocounter{oekaki}{-1}
\label{oe:fRijij}\caption{}\end{oekaki}

By an L-operator we mean a matrix
$L(u) = [L(u)_j^i]_{i,j=1,\cdots n}$
of operators (noncommutative letters) that satisfies

\begin{equation}
\check{R}(u-v) L(u) \otimes L(v)
= L(v) \otimes L(u) \check{R}(u-v),
\label{eq:RLL=LLR}
\end{equation}
where $\check{R}(u) := PR(u)$.

For Belavin's R-matrix we shall construct
such an L-operator in the following way.
Let ${\bf h}^*$ be the weight space for $sl_n(\cpx)$.
We realize ${\bf h}^*$
in $\cpx^n = \oplus_{i=1,\cdots, n}\cpx \epsilon_i$,
$<\epsilon_i,\epsilon_j>=\delta_{i,j}$, as the orthogonal complement to
$\sum_{i=1,\cdots,n} \epsilon_i$.
We denote the orthogonal projection of $\epsilon_i$ by $\bar{\epsilon}_i$.
For each $\lambda,\mu \in \mbox{\bf h}^*$ and $j=1,\cdots, n$ we
define the so-called intertwining vectors
(\cite{Bax73} for $n=2$ case and generalized in \cite{JMO88})
\begin{equation}
\ov     {u}
        {\mu}
        {\lambda}
        _j
:=
\left\{\begin{array}{cl}
        {\theta_j(\frac{u}{n}-<\lambda,\bar{\epsilon}_k>)}
         /{\sqrt{-1}\eta(\tau)}
        & :\mu-\lambda=\hbar\bar{\epsilon}_k
                \quad {\rm for \ some}\ k=1,\cdots,n,\\
        0
        & :{\rm otherwise}\\
\end{array}\right.
\label{eq:defov}
\end{equation}
where
$$
\theta_{j}(u):= \theta_{\frac{n}{2}-j,n}(u+\frac{1}{2},\tau)
=
        \sum_{\mu \in \frac{n}{2}-j+n{\bf Z}}^{}
        {\exp 2\pi i\left[\mu (u+\frac{1}{2}) + \frac{\mu^2}{2n}\tau
          \right]}
$$
and
$\eta(\tau):= p^{1/24}\prod_{m=1}^\infty (1-p^m)$
denotes the Dedekind eta function with $p:=\exp 2\pi i\tau$.


\begin{oekaki}\refstepcounter{oekaki}\addtocounter{oekaki}{-1}
\label{oe:foutitv}\caption{}\end{oekaki}
\noindent
Let further define
$
\iv
        {u}
        {\mu+\hbar \bar{\epsilon}_k,}
        {\mu}
        ^j
$
to be the entry in the inverse matrix to
$[
\ov{u}
        {\mu+\hbar \bar{\epsilon}_k}
        {\mu}_j
]
_{j,k=1,\cdots,n}
$,
namely \cite{QF}
\begin{equation}
\sum_{j=1}^n
{\bar{\phi}(u)}
        ^{\mu +\hbar\bar{\epsilon}_k}
        _{\mu}
        {}^j
\ov
        {u}
        {\mu +\hbar\bar{\epsilon}_{k'}}
        {\mu}
        _j
=
\delta_{k,k'}
,
\qquad
\sum_{k=1}^n
\ov
        {u}
        {\mu+\hbar\bar{\epsilon_k}}
        {\mu}
        _j
{\bar{\phi}(u)}
        ^{\mu+\hbar\bar{\epsilon_k}}
        _{\mu}
        {}^{j'}
=
\delta_{j,j'}
{}.
\label{eq:dualityofitv}
\end{equation}

The pictorial expression for $\phi, \bar{\phi}$ are given
in Figs.\ref{oe:foutitv},\ref{oe:fincitv} respectively,
and then (\ref{eq:dualityofitv}) can be written as in
Fig.\ref{fig:dualityofitv}.
\begin{oekaki}[h]\refstepcounter{oekaki}\addtocounter{oekaki}{-1}
\label{oe:fincitv}\caption{}\end{oekaki}

\begin{oekaki}[h]\refstepcounter{oekaki}\addtocounter{oekaki}{-1}
\label{fig:dualityofitv}\caption{}\end{oekaki}

\noindent
Then generalizing a result in the paper \cite{Sklyanin}(II),
we have

\begin{thm}
[\cite {has}, \cite{has93}]
For a function $f$ on ${\bf h}^*$, put
\begin{equation}
(L(c|u)^i_j f)(\mu):=
\sum_{k=1}^{n}
\ov     {u+c\hbar}
        {\mu+\hbar\bar{\epsilon}_k}
        {\mu}
        _j
\iv
        {u}
        {\mu+\hbar \bar{\epsilon}_k,}
        {\mu}
        ^i
f(\mu+\hbar\bar{\epsilon_k})
\label{eq:defofL}
\end{equation}
(Fig.\ref{fig:ffactorL}).
Then for any $c\in \cpx$, the collection of difference operators
$L(c|u)=[L(c|u)^i_j]_{i,j=1\cdots,n}$
satisfies the desired relation (\ref{eq:RLL=LLR}).
i.e., $L(c|u)$ gives a 1-parameter $(c)$ family of L-operators.
\label{thm:Lop}
\end{thm}

\begin{oekaki}\refstepcounter{oekaki}\addtocounter{oekaki}{-1}
\label{fig:ffactorL}\caption{}\end{oekaki}

{\bf Remark.}
This type of ``factorized'' L-operator first appeared in
 \cite{IK}.
It was for the trigonometric two by two R-matrix case, namely the
six vertex model corresponding to the algebra $U_q(\widehat{sl}_2)$.
Later, in connection with the
chiral-Potts and the Kashiwara-Miwa solutions
of the Yang-Baxter equation,
generalization is given in this form in the
$U_q(\widehat{sl}_n)$ case by \cite{BKMS} and the elliptic case by
\cite{QF} and \cite{has}.

\medskip


Let us give a brief account for the proof of this theorem.
In the proof, %
the so-called $A_{n-1}^{(1)}$ face model Boltzmann weight
\cite{JMO88}
(abbrev. face weight)
plays a central role.
The face weight is a collection of the quantities (Figure {fig:ffacewt})

$$
\pw
\left[\begin{array}{ccc}
                        & \mu&\\
                \lambda &u&\nu\\
                        & \mu'& \\\end{array}\right]
$$
\begin{oekaki}[h]\refstepcounter{oekaki}\addtocounter{oekaki}{-1}
\label{fig:ffacewt}\caption{}\end{oekaki}

\noindent
defined for
$\lambda, \mu, \mu', \nu \in {\bf h}^*$ and $u\in \cpx$.
Explicitly they are given by the formulas
$$
\pw
\left[\begin{array}{ccc}
                        & \lambda +\hbar\bar{\epsilon_{i}} & \\
                \lambda &u& \lambda + 2\hbar\bar{\epsilon_{i}} \\
                        & \lambda +\hbar\bar{\epsilon_{i}} &
                        \\\end{array}\right
]
:= \frac {\theta(u+{\hbar})}{\theta({\hbar})},
$$
$$
\pw
\left[\begin{array}{ccc}
                & \lambda +\hbar\bar{\epsilon_{i}} & \\
        \lambda
        &u&\lambda+\hbar(\bar{\epsilon_{i}}+\bar{\epsilon_{j}})\\
                & \lambda +\hbar\bar{\epsilon_{i}} &
                \\\end{array}\right]
:= \frac {\theta(-u+\lambda_{ij})}{\theta(\lambda_{ij})}
\; (i\neq j),
$$
$$
\pw
\left[\begin{array}{ccc}
                & \lambda +\hbar\bar{\epsilon_{i}} & \\
        \lambda &u&\lambda + \hbar(\bar{\epsilon_{i}}
        +\bar{\epsilon_{j}})\\
                & \lambda +\hbar\bar{\epsilon_{j}} &
                \\\end{array}\right]
:= \frac {\theta(u)}{\theta({\hbar})}
        \frac {\theta({\hbar}+\lambda_{ij})}{\theta(\lambda_{ij})}
\; (i\neq j),
$$
where
$
\lambda_{ij} := <\lambda, \bar{\epsilon_i} - \bar{\epsilon_j}>,
$
and
\begin{equation}
\theta(u):= \theta_{\frac{1}{2},1}(u+\frac{1}{2})
=
\sqrt{-1}p^{1/8}(z^{1/2}-z^{-1/2})
\prod_{m=1}^\infty
(1-zp^m)(1-z^{-1}p^m)(1-p^m)
\label{eq:tripleprod}
\end{equation}
$
(z^{1/2}=\exp \pi\sqrt{-1}u, p=\exp 2\pi\sqrt{-1}\tau)
$
denotes the Jacobi theta function.
For the other configurations of $\lambda,\mu,\mu'$ and $\nu$,
the face weight is set to 0:
$$
\pw
\left[\begin{array}{ccc}
                        & \mu&\\
                \lambda &u&\nu\\
                        & \mu'& \\\end{array}\right]
:= 0.
$$
The face weight and Belavin's R matrix are related by the following
relation (Fig.\ref{fig:fovprop}),
%
%
\begin{equation}
\sum_{i,j=1}^n
R(u-v)^{i,j}_{i',j'}
\ov     {u}
        {\mu}
        {\lambda}
        _{i}
\otimes
\ov     {v}
        {\nu}
        {\mu}
        _{j}
=
\sum_{\mu'}
\ov     {v}
        {\mu'}
        {\lambda}
        _{j'}
\otimes
\ov     {u}
        {\nu}
        {\mu'}
        _{i'}
\pw
\left[\begin{array}{ccc}
                        & \mu & \\
                \lambda & u-v &\nu\\
                        & \mu'& \\\end{array}\right]
\label{eq:ovintertwins}
\end{equation}
\begin{oekaki}\refstepcounter{oekaki}\addtocounter{oekaki}{-1}
\label{fig:fovprop}\caption{}\end{oekaki}
\noindent
i.e. they are ``intertwined'' by the intertwining
vectors.
It follows that the face weight satisfies the so-called face
Yang-Baxter equation (see (\ref{eq:faceYBE})).
By the duality relation (\ref{eq:dualityofitv}),
(\ref{eq:ovintertwins}) immediately implies
\begin{equation}
\sum_{i',j'=1}^n
\iv
        {v}
        {\mu}
        {\lambda}
        ^{j'}
\otimes
\iv
        {u}
        {\nu}
        {\mu}
        ^{i'}
R(u-v)^{i,j}_{i',j'}
=
\sum_{\mu'}
\pw
        \left[\begin{array}{ccc}
                        & \mu'& \\
                \lambda & u-v &\nu\\
                        & \mu & \\
        \end{array}\right]
\iv
        {u}
        {\mu'}
        {\lambda}
        ^{i}
\otimes
\iv
        {v}
        {\nu}
        {\mu'}
        ^{j}
\label{eq:ivintertwins}
\end{equation}
\begin{oekaki}\refstepcounter{oekaki}\addtocounter{oekaki}{-1}
\label{fig:fivprop}\caption{}\end{oekaki}

\noindent
(Fig. \ref{fig:fivprop})
which is in fact equivalent to (\ref{eq:ovintertwins}) itself.
Combining these two relations, we have

\begin{eqnarray*}
&&
\sum_{i',j';\mu}
\ov     {v+c\hbar}
        {\mu}
        {\lambda}
        _{j''}
\iv
        {v}
        {\mu}
        {\lambda}
        ^{j'}
\otimes
\ov     {u+c\hbar}
        {\nu}
        {\mu}
        _{i''}
\iv
        {u}
        {\nu}
        {\mu}
        ^{i'}
R(u-v)^{i,j}_{i',j'}
\\
&\stackrel{(\ref{eq:ivintertwins})}{=}&
\sum_{\mu,\mu'}
\ov     {v+c\hbar}
        {\mu}
        {\lambda}
        _{j''}
\otimes
\ov     {u+c\hbar}
        {\nu}
        {\mu}
        _{i''}
\pw
        \left[\begin{array}{ccc}
                        & \mu'& \\
                \lambda & u-v &\nu\\
                        & \mu & \\
        \end{array}\right]
\iv
        {u}
        {\mu'}
        {\lambda}
        ^{i}
\otimes
\iv
        {v}
        {\nu}
        {\mu'}
        ^{j}
\\
&\stackrel{(\ref{eq:ovintertwins})}{=}&
\sum_{i',j';\mu'}
R(u-v)^{i',j'}_{i'',j''}
\ov     {u+c\hbar}
        {\mu'}
        {\lambda}
        _{i'}
\iv
        {u}
        {\mu'}
        {\lambda}
        ^{i}
\otimes
\ov     {v+c\hbar}
        {\nu}
        {\mu'}
        _{j'}
\iv
        {v}
        {\nu}
        {\mu'}
        ^{j}
,
\end{eqnarray*}
(Fig. \ref{fig:fRLL_LLR}) and the theorem follows.
\hfill\fb
\begin{oekaki}\refstepcounter{oekaki}\addtocounter{oekaki}{-1}
\label{fig:fRLL_LLR}\caption{}\end{oekaki}

\section{Fusion Procedure}
\subsection{Overview}
Let $A(R)$ be the bialgebra generated by the formal letters
$\{L(u)^i_j|u\in \cpx, i,j=1,\cdots,n\}$
with respect to the relation
(\ref{eq:RLL=LLR})
and let ${\cal O}({\bf h}^*)$ be the ring of meromorphic functions
on ${\bf h}^*$.
Recall
$V=\oplus_{j=1,\cdots,n} \cpx e^j \simeq \cpx^n$.
Our L-operator (\ref{eq:defofL}),
giving a representation
$L(u)_j^i \mapsto L(c|u)_j^i$
of $A(R)$ by the definition,
is an endomorphism on the space
$
V \otimes {\cal O}({\bf h}^*)
$
:
$$
L(c|u) \in {\rm End}(V \otimes {\cal O}({\bf h}^*)).
$$
Here the first space $V = \cpx^n$  should be regarded as
the defining comodule (vector ``co''\-repre\-sen\-ta\-tion) for $A(R)$
: the comodule structure depends on the spectral parameter $u$
and is given by
$$
\fb_u :
V \ni e^i \mapsto \sum_{j=1}^n  e^j\otimes L(u)^i_j \in V\otimes A(R).
$$
We denote the space $V$ endowed with this comodule structure by
$V(\fb_u).$
We can consider more complicated comodules
for $A(R)$ as well.
That is, for each Young diagram $Y$ and a parameter $u$ we can construct
an $A(R)$-comodule $V(Y_u)$
whose dimension is just the same as for the $GL_n$-module that
corresponds to $Y$.
In contrast with the trigonometric R-matrix case,
where we have the quantized enveloping algebra
$U_\hbar(A_{n-1}^{(1)})$
and its universal R-matrix
as the origin of the tensor category structure,
in the elliptic R-matrix case such an algebraized theory
does not seem to be available until now.
Yet we can utilize a constructive definition of  $V(Y_u)$
known as the fusion technique
\cite{KRS}\cite{C}.
This is done by taking appropriate sub/quotient of the
tensor comodule
$
V(\fb_{u_1}\otimes\cdots\otimes\fb_{u_k})
=
V(\fb_{u_1})\otimes\cdots\otimes V(\fb_{u_k})
\simeq V^{\otimes k}$.
Applying this fusion technique to our case,
we get a collection of difference operators
$$
L(c|Y_u) \in {\rm End}(V(Y_u) \otimes {\cal O}({\bf h}^*))
$$
for each $Y$, as well as the ``fused R-matrices''
$
\check{R}({Y_u,Y'_v}):
V(Y_u)\otimes V(Y'_v) \rightarrow V(Y'_v) \otimes V(Y_u).
$
The latter gives an isomorphism between the two
$A(R)$-comodules
$V(Y_u)\otimes V(Y'_v)$ and $V(Y'_v) \otimes V(Y_u)$
and it follows that we have the relation
\begin{equation}
\check{R}({Y_u,Y'_v})\, L(c|Y_u) \otimes L(c|Y'_v)
=
L(c|Y'_v) \otimes  L(c|Y_u)\, \check{R}({Y_u,Y'_v})
{}.
\label{eq:RLL'=L'LR}
\end{equation}

\subsection{Antisymmetric case}
For our purpose of constructing Macdonald type commuting difference operators,
important are the fused L operators corresponding to the Young diagrams of
vertical boxes that we are now going to describe.

\medskip
Let us denote the R matrix
$\check{R}(u-v)$
(\ref{eq:RTformula}) considered as an operator
$V(\fb_u) \otimes V(\fb_v) \rightarrow V(\fb_v)\otimes V(\fb_u)$
by
$\PR{\fb_u,\fb_v}.$
The Yang-Baxter equation reads
\begin{eqnarray*}
(1 \otimes \PR{\fb_{u},\fb_{v}})
(\PR{\fb_{u},\fb_{w}} \otimes 1)
(1 \otimes \PR{\fb_{v},\fb_{w}})
&=&
(\PR{\fb_{v},\fb_{w}} \otimes 1)
(1 \otimes \PR{\fb_{u},\fb_{w}})
(\PR{\fb_{u},\fb_{v}} \otimes 1)
\\
:
V(\fb_u) \otimes V(\fb_v) \otimes V(\fb_w)
&\rightarrow&
V(\fb_w) \otimes V(\fb_v)\otimes V(\fb_u).
\label{eq:vertexYBE}
\end{eqnarray*}
We put (Fig. \ref{fig:ffusedR1})
\begin{eqnarray}
\lefteqn{
\PR
{\fb_{u_1}\otimes\fb_{u_2}\otimes\cdots\otimes\fb_{u_{k}},\fb_v}
:=
(\PR{\fb_{u_1},\fb_{v}}
                        )^{1,2}
(\PR{\fb_{u_2},\fb_{v}}
                        )^{2,3}
\cdots
(\PR{\fb_{u_k},\fb_{v}}
                        )^{k,k+1}
}
\label{eq:fusedR_1}
\\
&:&
V(\fb_{u_1}\otimes\fb_{u_2}\otimes\cdots\otimes\fb_{u_{k}})
\otimes V(\fb_v)
\rightarrow
 V(\fb_v)\otimes
V(\fb_{u_1}\otimes\fb_{u_2}\otimes\cdots\otimes\fb_{u_{k}}).
\nonumber
\end{eqnarray}
\begin{oekaki}\refstepcounter{oekaki}\addtocounter{oekaki}{-1}
\label{fig:ffusedR1}\caption{}\end{oekaki}

\noindent
The superscripts show the tensor components on which each operator
acts. The resulting operator (\ref{eq:fusedR_1})
is the operator which sends the rightmost tensor component to the left.
We also put (Fig. \ref{fig:ffusedR2})
\begin{eqnarray}
\lefteqn{
\PR
{\fb_{u_1}\otimes\cdots\otimes\fb_{u_{k}}
,\fb_{v_1}\otimes\cdots\otimes\fb_{v_{l}}}
}
\nonumber\\
&:=&
(\PR
{\fb_{u_1}\otimes\cdots\otimes\fb_{u_{k}},\fb_{v_l}}
                                                )^{l\cdots k+l-1;k+l}
\cdots
\nonumber\\
&&
(\PR
{\fb_{u_1}\otimes\cdots\otimes\fb_{u_{k}},\fb_{v_2}}
                                                )^{2\cdots k+1;k+2}
(\PR
{\fb_{u_1}\otimes\cdots\otimes\fb_{u_{k}},\fb_{v_1}}
                                                )^{1\cdots k;k+1}
\label{eq:composedR}
\\
&:&
V(\fb_{u_1}\otimes\cdots\otimes\fb_{u_{l}})
\otimes
V(\fb_{v_1}\otimes\cdots\otimes\fb_{v_{k}})
\rightarrow
V(\fb_{v_1}\otimes\cdots\otimes\fb_{v_{k}})
\otimes
V(\fb_{u_1}\otimes\cdots\otimes\fb_{u_{l}})
{}.
\nonumber
\end{eqnarray}
\begin{oekaki}\refstepcounter{oekaki}\addtocounter{oekaki}{-1}
\label{fig:ffusedR2}\caption{}\end{oekaki}

\noindent
For $k=1,\cdots,n$,
let $1^k$ be the Young diagram of vertical $k$ boxes \cite{Macbook}
(
$1^1=\fb$
%
).
Then the
fusion operator by Cherednik
associated with $1^k$
is given by the comodule map corresponding to ``the half twist''
(Fig. \ref{fig:fvproj})
\begin{eqnarray}
\pi_{1^k}
&:=&
(\PR
{\fb_{u_1},\fb_{u_2}}
                        )^{k-1;k}
\cdots
\nonumber\\
&&
(\PR
{\fb_{u_1}\otimes\fb_{u_2}\otimes\cdots\otimes\fb_{u_{k-2}},\fb_{u_{k-1}}}
)       ^{2\cdots k-1;k}
(\PR
{\fb_{u_1}\otimes\fb_{u_2}\otimes\cdots\otimes\fb_{u_{k-1}},\fb_{u_k}}
)       ^{1\cdots k-1;k}
\nonumber
\\
&:&
V(\fb_{u-(k-1){\hbar}}\otimes\cdots\otimes \fb_{u})
\rightarrow
V(\fb_{u}\otimes\cdots\otimes \fb_{u-(k-1){\hbar}})
\label{eq:defpi}
\end{eqnarray}
where the spectral parameters are specialized as
\begin{equation}
(u_1,\cdots,u_k)
=
(u-(k-1){\hbar},\cdots,u-{\hbar},u)
\label{eq:specialparam}
\end{equation}
\begin{oekaki}\refstepcounter{oekaki}\addtocounter{oekaki}{-1}
\label{fig:fvproj}\caption{}\end{oekaki}

\noindent
so that the rank of the operator $\pi_{1^k}$ degenerates.
Since the operator $\PR{\fb_u,\fb_v}$ actually depends only in the
difference $u-v$,
the operator
$\pi_{1^k}$
does not depend on $u$.
The reason why this specialization
(\ref{eq:specialparam})
of parameters gives rise to the degeneration is as follows.
First of all, from (\ref{eq:RTformula}) we have
\begin{equation}
\check{R}(-\hbar)_{i'j'}^{ij}=-\check{R}(-\hbar)_{j'i'}^{ij},
\label{eq:imR=ext2}
\end{equation}
showing ${\rm Im}(\check{R}(-\hbar))\subset \wedge^2(\cpx^n)$.
Actually for a generic value of $\hbar$ one can show
$
{\rm Im}(\check{R}(-\hbar)) = \wedge^2(\cpx^n)
$.
Similarly, in $\pi_{1^k}$ one finds a degenerated factor
$\check{R}(-\hbar)$ for each adjacent pair of tensor factors.
With some braid manipulation, this implies
the antisymmetricity of the matrix element of  $\pi_{1^k}$
with respect to the permutation of the outgoing indices.
This shows
$
{\rm Im} \pi_{1^k} \subset \wedge^k(\cpx^n)
$
and in fact
${\rm Im} \pi_{1^k} \simeq \wedge^k(\cpx^n)$
for generic $\hbar$.
We put
$$
V({1^k}_u):=
\pi_{1^k}(V(\fb_{u-(k-1){\hbar}}\otimes\cdots\otimes\fb_{u}))
\subset
V(\fb_{u}\otimes\cdots\otimes\fb_{u-(k-1){\hbar}}).
$$

\medskip
The fused R matrix as the comodule map
$
V(1^k_u)\otimes V(1^{k'}_v) \rightarrow V(1^{k'}_v)\otimes V(1^k_u)
$
is now defined as the restriction of the composed R matrix
(\ref{eq:composedR}),
$$
\PR
{1^k_u,1^{k'}_v}
:=
\left.
\PR
{\fb_{u}\otimes\cdots\otimes\fb_{u-(k-1){\hbar}}
,\fb_{v}\otimes\cdots\otimes\fb_{v-(k'-1){\hbar}}}
\right|
_{V(1^k_u)\otimes V(1^{k'}_v)}.
$$
\medskip
We also define the fused L operator
$
L(1^k_u)^I_{I'}\in A(R)
$
to be the matrix element of the comodule structure map.
That is, choosing a basis $\{e^I\}_I \subset V(1^k_u)$, we write
$$
1^k_u :
 V(1^k_u) \ni e^I
 \mapsto
 \sum_{I'}  e^{I'}\otimes L(1^k_u)^I_{I'}\in V(1^k_u)\otimes A(R) ,
$$
or equivalently, if we denote the dual basis to
$\{e^I\}\subset V(1^k_u)$
by
$\{e_I\}\subset V(1^k_u)^*$,
we put
$$
L(1^k_u)^I_{I'}:= <e_{I'}, 1^k_u(e^I)>.
$$
Since $V(1^k_u) \simeq \wedge^k V$ as a vector space,
we can choose the basis elements in $V(1^k_u)$ as
$$e^I
:=
\sum_{\sigma\in S(k)}{\rm sgn}(\sigma)
e^{i_{\sigma(1)}}\otimes \cdots\otimes e^{i_{\sigma(k)}}
\,\,
(I=\{i_1,\cdots,i_k|i_1<\cdots<i_k\} \subset \{1,\cdots,n\}),
$$
where $S(k)$ denotes the symmetric group of order $k$.
Let $\{e_i\}\subset V^*$ be the dual basis for $\{e^i\}\subset V$,
then for the above basis
$\{e^I\}\subset V(1^k_u)$
the dual basis is just given by
$$
e_I := e_{i_1}\otimes\cdots\otimes e_{i_k}|_{V(1^k_u)}
\,\,
(I=\{i_1,\cdots,i_k|i_1<\cdots<i_k\} \subset \{1,\cdots,n\}).
$$
With respect to this choice of basis, the fused L operator is
given by
\begin{eqnarray}
L(1^k_u)^I_{I'}
&=& <e_{I'}, 1^k_u(e^I)>
\nonumber
\\
&=&
<
e_{i'_1}\otimes\cdots\otimes e_{i'_k},
1^k_u
(\sum_{\sigma\in S(k)}{\rm sgn}(\sigma)
  e^{i_{\sigma(1)}}\otimes \cdots\otimes e^{i_{\sigma(k)}}
)
>
\nonumber
\\
&=&
<
e_{i'_1}\otimes\cdots\otimes e_{i'_k},
\fb_{u}\otimes\cdots\otimes\fb_{u-(k-1){\hbar}}
(\sum_{\sigma\in S(k)}{\rm sgn}(\sigma)
e^{i_{\sigma(1)}}\otimes \cdots\otimes e^{i_{\sigma(k)}}
)
>
\nonumber
\\
&=&
\sum_{\sigma\in S(k)}{\rm sgn}(\sigma)
<e_{i'_1},\fb_{u}(e^{i_{\sigma(1)}})>
\cdots
<e_{i'_k},\fb_{u-(k-1){\hbar}}(e^{i_{\sigma(k)}})>
\nonumber
\\
&=&
\sum_{\sigma\in S(k)}{\rm sgn}(\sigma)
L(u)_{i'_1}^{i_{\sigma(1)}}
\cdots
L(u-(k-1){\hbar})_{i'_k}^{i_{\sigma(k)}}
\label{eq:fusedL}
\end{eqnarray}
where $I=\{i_1<\cdots<i_k\}$ and $I'=\{i'_1<\cdots<i'_k\}.$
Applying the representation (\ref{eq:defofL}) to this element
we obtain a formula for the fused L operator
acting on the space of functions on ${\bf h}^*$,
\begin{equation}
A(R)\ni L(1^k_u)^I_{I'}
\mapsto
 L(c|1^k_u)^I_{I'}
:=
\sum_{\sigma\in S(k)}{\rm sgn}(\sigma)
L(c|u)_{i'_1}^{i_{\sigma(1)}}
\cdots
L(c|u-(k-1){\hbar})_{i'_k}^{i_{\sigma(k)}}
\in
{\rm End}
 ({\cal O}({\bf h}^*)).
\label{eq:fuseddiffL}
\end{equation}
 From the construction, they enjoy the relation (\ref{eq:RLL'=L'LR}).

\subsection{Fused L operators from fused intertwining vectors}
\medskip
The basic L operator (\ref{eq:defofL}) is defined
as the composition of two quantities, the intertwining vectors.
The fused L operator (\ref{eq:fuseddiffL})
also have this factorized nature as we explain below.

\medskip
First we introduce ``the space of paths'' and their fusions
to formulate the face weights as linear maps on appropriate vector
spaces.
We put
$$
\path     {\fb}
        {\mu}
        {\lambda}
:\cong
\left\{\begin{array}{cl}
        {\cpx}  &
                  :\mu =\lambda + \hbar\bar{\epsilon_i} {\rm \ for \ some
                    }\ i,\\
        0       & :{\rm otherwise.}\\
\end{array}\right.
$$
(Fig.\ref{fig:f1steppath})
and denote by
$
e_{\lambda}^{\mu}
$
the basis of the one dimensional space
$
\path     {\fb}
        {\mu}
        {\lambda}
$
when
$
\mu =\lambda + \hbar\bar{\epsilon_i}
$
for some $i$,
and otherwise we set
$
e_{\lambda}^{\mu}=0.
$
\begin{oekaki}\refstepcounter{oekaki}\addtocounter{oekaki}{-1}
\label{fig:f1steppath}\caption{}\end{oekaki}

\noindent
For each $u \in \cpx$ we consider the copy
$
\path     {\fb_u}
        {\mu}
        {\lambda}
$
of
$
\path     {\fb}
        {\mu}
        {\lambda}
$
and define (Fig. \ref{fig:flongpath})
$$
\path     {\fb_{u_1}\otimes\cdots\otimes\fb_{u_k}}
        {\nu}
        {\lambda}
:=
\bigoplus_{\mu_1,\cdots,\mu_{k-1}}
\path     {\fb_{u_1}}
        {\mu_1}
        {\lambda}
\otimes
\path     {\fb_{u_2}}
        {\mu_2}
        {\mu_1}
\otimes \cdots \otimes
\path     {\fb_{u_{k-1}}}
        {\nu}
        {\mu_{k-1}}
,
$$
$$
\path     {\fb_{u_1}\otimes\cdots\otimes\fb_{u_k}}
        {}{}
:=
\bigoplus_
        {\lambda,\nu}
\path     {\fb_{u_1}\otimes\cdots\otimes\fb_{u_k}}
        {\nu}
        {\lambda}
{}.
$$
\begin{oekaki}\refstepcounter{oekaki}\addtocounter{oekaki}{-1}
\label{fig:flongpath}\caption{}\end{oekaki}

\noindent
These are the space of ``admissible paths'' in \cite{An1face},
in which the face weight and their fusions are formulated as
linear operators on these spaces as follows.
For
$
\og{e}  {}
        {\lambda}
        {\mu}
\in
\path   {\fb_{u}}
        {\mu}
        {\lambda}
$
and
$
\og{e}  {}
        {\mu}
        {\nu}
\in
\path   {\fb_{v}}
        {\nu}
        {\mu}
$
we put
$$
\PW{\fb_u,\fb_v}
(e_{\lambda}^{\mu} \otimes e_{\mu}^{\nu})
:= \sum_{\mu'}
        e_{\lambda}^{\mu'} \otimes e_{\mu'}^{\nu}
\pw
\left[\begin{array}{ccc}
                        & \mu   & \\
                \lambda & u-v   &\nu\\
                        & \mu'  & \\\end{array}\right],
$$
thereby define the face operator (Fig.\ref{fig:ffaceop})
$
\PW{\fb_u,\fb_v}
:
\path   {\fb_{u}\otimes\fb_{v}}
        {}{}
\rightarrow
\path   {\fb_{v}\otimes\fb_{u}}
        {}{}
,\quad
\path   {\fb_{u}\otimes\fb_{v}}
        {\nu}
        {\lambda}
\rightarrow
\path   {\fb_{v}\otimes\fb_{u}}
        {\nu}
        {\lambda}
{}.
$
\begin{oekaki}\refstepcounter{oekaki}\addtocounter{oekaki}{-1}
\label{fig:ffaceop}\caption{}\end{oekaki}

\noindent
With these definitions,
the Yang-Baxter equation of face type reads as follows
(Fig.\ref{fig:ffaceYBE}):
\begin{eqnarray}
(1 \otimes \PW{\fb_{u},\fb_{v}})
(\PW{\fb_{u},\fb_{w}} \otimes 1)
(1 \otimes \PW{\fb_{v},\fb_{w}})
&=&
(\PW{\fb_{v},\fb_{w}} \otimes 1)
(1 \otimes \PW{\fb_{u},\fb_{w}})
(\PW{\fb_{u},\fb_{v}} \otimes 1)
\nonumber
\\
:
\path   {\fb_{u}\otimes\fb_{v}\otimes\fb_{w}}
        {\nu}
        {\lambda}
&\rightarrow&
\path   {\fb_{w}\otimes\fb_{v}\otimes\fb_{u}}
        {\nu}
        {\lambda}
{}.
\label{eq:faceYBE}
\end{eqnarray}
\begin{oekaki}\refstepcounter{oekaki}\addtocounter{oekaki}{-1}
\label{fig:ffaceYBE}\caption{}\end{oekaki}

\noindent
The fusion procedure works for the face weight as well.
We put
\begin{eqnarray}
\lefteqn{
\PW
{\fb_{u_1}\otimes\fb_{u_2}\otimes\cdots\otimes\fb_{u_{k}},\fb_v}
:=
(\PW{\fb_{u_1},\fb_{v}}
                        )^{1,2}
(\PW{\fb_{u_2},\fb_{v}}
                        )^{2,3}
\cdots
(\PW{\fb_{u_k},\fb_{v}}
                        )^{k,k+1}
}
\nonumber\\
&:&
{\cal P}(\fb_{u_1}\otimes\fb_{u_2}\otimes\cdots\otimes\fb_{u_{k}}
\otimes \fb_v)
\rightarrow
{\cal P}(\fb_v\otimes
\fb_{u_1}\otimes\fb_{u_2}\otimes\cdots\otimes\fb_{u_{k}})
,
\label{eq:fusedW_1}
\end{eqnarray}
\begin{eqnarray}
\lefteqn{
\PW
{\fb_{u_1}\otimes\cdots\otimes\fb_{u_{k}}
,\fb_{v_1}\otimes\cdots\otimes\fb_{v_{l}}}
}
\nonumber\\
&:=&
(\PW
{\fb_{u_1}\otimes\cdots\otimes\fb_{u_{k}},\fb_{v_l}}
                                                )^{k\cdots k+l-1;k+l}
\cdots
\nonumber\\
&&
\PW
{\fb_{u_1}\otimes\cdots\otimes\fb_{u_{k}},\fb_{v_2}}
                                                )^{2\cdots k+1;k+2}
(\PW
{\fb_{u_1}\otimes\cdots\otimes\fb_{u_{k}},\fb_{v_1}}
                                                )^{1\cdots k;k+1}
\label{eq:composedW}
\\
&:&
{\cal P}
(\fb_{u_1}\otimes\cdots\otimes\fb_{u_{l}}
\otimes
\fb_{v_1}\otimes\cdots\otimes\fb_{v_{k}})
\rightarrow
{\cal P}
(\fb_{v_1}\otimes\cdots\otimes\fb_{v_{k}}
\otimes
\fb_{u_1}\otimes\cdots\otimes\fb_{u_{l}})
,
\nonumber
\end{eqnarray}
and the fusion operator for $1^k$ is given by (Fig.\ref{fig:ffproj})
\begin{eqnarray}
\Pi_{1^k}
&:=&
(\PW
{\fb_{u_1},\fb_{u_2}}
                        )^{k-1;k}
\cdots
\nonumber\\
&&
(\PW
{\fb_{u_1}\otimes\fb_{u_2}\otimes\cdots\otimes\fb_{u_{k-2}},\fb_{u_{k-1}}}
)       ^{2\cdots k-1;k}
(\PW
{\fb_{u_1}\otimes\fb_{u_2}\otimes\cdots\otimes\fb_{u_{k-1}},\fb_{u_k}}
)       ^{1\cdots k-1;k}
\nonumber
\\
&&
\,:\,
{\cal P}(\fb_{u-(k-1){\hbar}}\otimes\cdots\otimes \fb_{u})
\rightarrow
{\cal P}(\fb_{u}\otimes\cdots\otimes \fb_{u-(k-1){\hbar}})
\label{eq:defPi}
\end{eqnarray}
with the same specialization (\ref{eq:specialparam}) of
$
(u_1,\cdots,u_k).
$
\begin{oekaki}\refstepcounter{oekaki}\addtocounter{oekaki}{-1}
\label{fig:ffproj}\caption{}\end{oekaki}

\noindent
Like $\pi_{1^k}$, $\Pi_{1^k}$ does not depend on $u$.
We denote the image of this operator by
$
\path{1^k}{}{}
:=\Pi_{1^k}
(
{\cal P}(\fb_{u-(k-1){\hbar}}\otimes\cdots\otimes \fb_{u})
)
\subset
{\cal P}(\fb_{u}\otimes\cdots\otimes \fb_{u-(k-1){\hbar}}).
$
For $J=\{j_1<\cdots<j_k\}$ put
$\bar{\epsilon}_J:=\bar{\epsilon_{j_1}}+\cdots+\bar{\epsilon_{j_d}}$.
We have
\begin{equation}
\path{1^k}{}{}
=
\bigoplus_{\lambda,J\subset\{1,\cdots,n\},|J|=k}
\path
        {1^k}
        {\lambda+\hbar\bar{\epsilon}_J}
        {\lambda}
,
\,
{\rm where}
\,
\path
        {1^k}
        {\lambda+\hbar\bar{\epsilon}_J}
        {\lambda}
:=
\path
        {\fb_{u}\otimes\cdots\otimes \fb_{u-(k-1){\hbar}}}
        {\lambda+\hbar\bar{\epsilon}_J}
        {\lambda}
\cap
\path{1^k}{}{}.
\label{eq:sumofpath}
\end{equation}
The property similar to (\ref{eq:imR=ext2})
for the face weight with specialized
spectral parameter $-\hbar$ implies that each subspace
$
\path
        {1^k}
        {\lambda+\hbar\bar{\epsilon}_J}
        {\lambda}
$
$
(J\subset\{1,\cdots,n\},|J|=k)
$
is one dimensional and spanned by the antisymmetric tensor
\begin{equation}
e(1^k)
       ^{\lambda+\hbar\bar{\epsilon}_J}
       _{\lambda}
:=
\sum_{\sigma\in S(k)}
{\rm sgn}(\sigma)
e_\lambda(j_{\sigma(1)},\cdots,j_{\sigma(k)})
\in
\path{1^k}
        {\lambda+\hbar\bar{\epsilon}_J}
        {\lambda}
\label{eq:fusedpath}
\end{equation}
where 
$$
e_\lambda(i_1,\cdots,i_k)
:=
e_\lambda^{\lambda+\hbar\bar{\epsilon}_{i_1}}
\otimes
e_{\lambda+\hbar\bar{\epsilon}_{i_1}}
 ^{\lambda+\hbar\bar{\epsilon}_{i_1}+\hbar\bar{\epsilon}_{i_2}}
\otimes\cdots\otimes
e_{\lambda+\hbar\bar{\epsilon}_{i_1}+\cdots+\hbar\bar{\epsilon}_{i_{k-1}}}
 ^{\lambda+\hbar\bar{\epsilon}_{i_1}+\cdots+\hbar\bar{\epsilon}_{i_k}}
{}.
$$
If we denote the dual basis for
$
e_\lambda^{\lambda+\hbar\bar{\epsilon}_{i}}
\in
\path{\fb_u}
        {\lambda+\hbar\bar{\epsilon}_i}
        {\lambda}
$
by
$
{e^*}_\lambda^{\lambda+\hbar\bar{\epsilon}_{i}}
\in
\left(
\path{\fb_u}
        {\lambda+\hbar\bar{\epsilon}_i}
        {\lambda}
\right)^*
$,
the dual basis for
$
e(1^k)
       ^{\lambda+\hbar\bar{\epsilon}_J}
       _{\lambda}
$
is given by
\begin{equation}
{e^*(1^k)}
       ^{\lambda+\hbar\bar{\epsilon}_J}
       _{\lambda}
\equiv
e^*_\lambda(j_1,\cdots,j_k)
:=
{e^*}_\lambda^{\lambda+\hbar\bar{\epsilon}_{j_1}}
\otimes
{e^*}_{\lambda+\hbar\bar{\epsilon}_{j_1}}
     ^{\lambda+\hbar\bar{\epsilon}_{j_1}+\hbar\bar{\epsilon}_{j_2}}
\otimes\cdots\otimes
{e^*}_{\lambda+\hbar\bar{\epsilon}_{j_1}+\cdots+\hbar\bar{\epsilon}_{j_{k-1}}}
     ^{\lambda+\hbar\bar{\epsilon}_{j_1}+\cdots+\hbar\bar{\epsilon}_{j_k}}
|_
{\path{1^k}
        {\lambda+\hbar\bar{\epsilon}_J}
        {\lambda}
}.
\label{eq:fusedpath*}
\end{equation}

\medskip
To define fused intertwining vectors,
we formulate the intertwining vectors as linear maps as follows.

\medskip
Outgoing intertwining vectors (Fig.\ref{fig:fov}):
$$
\ov     {\fb_{u}}
{}{}
:
\path
        {\fb_{u}}
{}{}
\rightarrow
        V(\fb_u)
,
\qquad
\ov     {\fb_{u}}
{}{}
\og{e}
        {}
        {\lambda}
        {\mu}
:=
\sum_j e^j
\,\cdot\,
 \ov   {\fb_{u}}
        {\mu}
        {\lambda}
        _j
{}.
$$
\begin{oekaki}\refstepcounter{oekaki}\addtocounter{oekaki}{-1}
\label{fig:fov}\caption{}\end{oekaki}

\noindent

Incoming intertwining vectors (Fig.\ref{fig:fiv}):
$$
\iv     {\fb_{u}}
{}{}
:
V(\fb_u)
\rightarrow
\path
        {\fb_{u}}
{}{}
,
\qquad
\iv     {\fb_{u}}
{}{}
e^j
:=
\og{e}
        {}
        {\lambda}
        {\mu}
\,\cdot\,
\iv   {\fb_{u}}
        {\mu}
        {\lambda}
        ^j
{}.
$$
\begin{oekaki}\refstepcounter{oekaki}\addtocounter{oekaki}{-1}
\label{fig:fiv}\caption{}\end{oekaki}
\noindent

\noindent
The intertwining property for these vectors read
respectively as follows.
\begin{equation}
\PR{\fb_u,\fb_v}
\,
\ov{\fb_u}{}{} \otimes \ov{\fb_v}{}{}
=
\ov{\fb_v}{}{} \otimes \ov{\fb_u}{}{}
\PW{\fb_u,\fb_v}
\,:\,
\path{\fb_u\otimes\fb_v}{}{}
\rightarrow
V(\fb_v\otimes\fb_u),
\label{eq:ov_as_op}
\end{equation}
\begin{equation}
\iv{\fb_v}{}{} \otimes \iv{\fb_u}{}{}
\PR{\fb_u,\fb_v}
=
\PW{\fb_u,\fb_v}
\iv{\fb_u}{}{} \otimes \iv{\fb_v}{}{}
\,:\,
V(\fb_u\otimes\fb_v)
\rightarrow
\path{\fb_v\otimes\fb_u}{}{}.
\label{eq:iv_as_eq}
\end{equation}
Putting
$$
\ov{\fb_{u_1}\otimes\cdots\otimes\fb_{u_k}}{}{}
:=\ov{\fb_{u_1}}{}{}\otimes\cdots\otimes\ov{\fb_{u_k}}{}{}
,\quad
\iv{\fb_{u_1}\otimes\cdots\otimes\fb_{u_k}}{}{}
:=\iv{\fb_{u_1}}{}{}\otimes\cdots\otimes\iv{\fb_{u_k}}{}{}
,
$$
the relations (\ref{eq:ov_as_op}), (\ref{eq:iv_as_eq})
respectively imply
\begin{equation}
\pi_{1^k}
\ov{\fb_{u_1}\otimes\cdots\otimes\fb_{u_k}}{}{}
=
\ov{\fb_{u_k}\otimes\cdots\otimes\fb_{u_1}}{}{}
\Pi_{1^k}
,
\quad
\iv{\fb_{u_k}\otimes\cdots\otimes\fb_{u_1}}{}{}
\pi_{1^k}
=
\Pi_{1^k}
\iv{\fb_{u_1}\otimes\cdots\otimes\fb_{u_k}}{}{}
\label{eq:piphi=phiPi}
\end{equation}
\begin{oekaki}\refstepcounter{oekaki}\addtocounter{oekaki}{-1}
\label{fig:piphi=phiPi}\caption{}\end{oekaki}
\noindent
(Fig.\ref{fig:piphi=phiPi}), where
$
(u_1,\cdots,u_k)
$
is as in
$
(\ref{eq:specialparam}).
$
Now we can define 
the fused intertwining vectors (outgoing and incoming) as
the maps
$$
\ov{1^k_u}{}{}
:=
\ov{\fb_{u_k}\otimes\cdots\otimes\fb_{u_1}}{}{}
|_{{\cal P}(1^k_u)}
:
{\cal P}(1^k_u) \rightarrow V(1^k_u)
,
$$
$$
\iv{1^k_u}{}{}
:=
\iv{\fb_{u_k}\otimes\cdots\otimes\fb_{u_1}}{}{}
|_{V(1^k_u)}
: V(1^k_u) \rightarrow {\cal P}(1^k_u).
$$
We also put
$$
\ov{1^k_u}{\mu}{\lambda}
:=
\ov{1^k_u}{}{}|_{\path{1^k_u}{\mu}{\lambda}}
,
\quad
\iv{1^k_u}{\mu}{\lambda}
:=
proj_{\lambda}^{\mu}\,\iv{1^k_u}{}{}
$$
where
$
proj_{\lambda}^{\mu}
:
\path{1^k}{}{}
\rightarrow
\path
        {1^k}
        {\nu}         
        {\lambda}
$
is the projection
onto the direct summand in (\ref{eq:sumofpath}).
(
$
proj_{\lambda}^{\mu}=0
$
unless
$
\mu={\lambda+\hbar\bar{\epsilon}_J}
$
for some $J$.
)

The factorized nature for our fused L operator
(\ref{eq:fuseddiffL})
is now stated as follows,

\begin{prop}
For $J=\{j_1<\cdots<j_k\}$, let
$\bar{\epsilon}_J:=\bar{\epsilon_{j_1}}+\cdots+\bar{\epsilon_{j_d}}$
and
$
T_J^\hbar f(\lambda):= f(\lambda+\hbar\bar{\epsilon}_J).
$
We have
\begin{equation}
\left( L(c|1^k_u)^I_{I'} f\right) (\lambda)
=
 \sum_{J\subset\{1,\cdots,n\},|J|=k}
  \ov   {1^k_{u+c\hbar}}
        {\lambda+\hbar\bar{\epsilon}_J}
        {\lambda}
        _{I'}
\circ
  \iv
        {1^k_u}
        {\lambda+\hbar\bar{\epsilon}_J}
        {\lambda}
        ^{,I}
 (T_J^\hbar f)(\lambda).
\label{eq:fusedLviaitv}
\end{equation}
Here $\circ$ stands for the composition of 
$
\iv        {1^k_u}{}{}
:
V(1^k_u) \rightarrow {\cal P}(1^k_u)
$
and
$
\ov        {1^k_{u+c\hbar}}{}{}
:{\cal P}(1^k_{u+c\hbar})
 \rightarrow V(1^k_u)
$
via the obvious identification
${\cal P}(1^k_u)\simeq {\cal P}(1^k_{u+c\hbar})$
,
super/subscripts $I,I'$ denote their matrix elements.

\label{prop:fusedLviaitv}
\end{prop}

Proof of this proposition is rather straightforward from the
definition and omitted.\hfill \fb

\section{The commuting family}

Now we are in a position to ask the question in
Section \ref{sec:intro}:
``what kind of operators arise as the traces of these
$L(c|Y_u)$ 's?".
We shall consider the case $Y=1^d$, the vertical $d$ boxes case.
Then
$L(c|1^d_u)$
is a matrix of size
${\rm dim} \wedge^d \cpx^n$
whose matrix elements are difference operators.

\begin{thm}[Main Theorem]

Let $M_{d}(c|u):=\tr_{1^d} L(c|1^d_u), \quad d=1,\cdots,n$.
We have
$$
M_{d}(c|u)
 =
\frac{\theta(u+\frac{dc\hbar}{n})}{\theta(u)}
\sum_
     {I\subset\{1,\cdots,n\},|I|=d}
\left(
 \prod
      _{s \notin I, t \in I}
  \frac
    {\theta(<\lambda,\bar{\epsilon}_s-\bar{\epsilon}_t> + \frac{c\hbar}{n})}
    {\theta(<\lambda,\bar{\epsilon}_s-\bar{\epsilon}_t>)}
\right)
  T_I^\hbar,
$$
where $\theta(u)$ is the Jacobi theta function (\ref{eq:tripleprod})
and
$T_I^\hbar$ stands for the $\hbar$-shift operator:
$$
(T_i^\hbar f)(\lambda):= f(\lambda +\hbar\bar{\epsilon}_i),
\quad
T_I^\hbar :=\prod_{i\in I} T_i^\hbar.
$$
These operators form a commuting difference system:
$$
[\,M_d(c|u), M_{d'}(c|v)\,]=0.
$$
\label{thm:trL_d}
\end{thm}

{\bf Remark.} It is an important point that these operators obviously commute,
as we mentioned in the introduction.
This is because the extended ``RLL=LLR" relation (\ref{eq:RLL'=L'LR})
can be rewritten as
$$
\check{R}({Y_u,Y'_v}) L(c|Y_u) \otimes L(c|Y'_v)
\check{R}({Y_u,Y'_v})^{-1}
=
L(c|Y'_v) \otimes  L(c|Y_u),
$$
and then taking the trace simply gives
$$
M(c|Y_u) M(c|Y'_v) = M(c|Y'_v) M(c|Y_u)
$$
where
$M(c|Y_u) := \tr_{V(Y)} L(c|Y_u)$.
This simple argument and the resulting operators,
``the commuting transfer matriecs'',
was effectively used in Baxter's analysis of the spin chain models
\cite{Bax71}. 
Thus our result can be said as an ideology:
\begin{center}
commuting transfer matrices = commuting difference system.
\end{center}

\medskip
{\bf Proof of Theorem \ref{thm:trL_d}.}
 From Proposition \ref{prop:fusedLviaitv},
the operator $M_d$ is of the form
\begin{eqnarray}
M_d(c|u)f(\lambda)
&=&
\sum_{I,J\subset\{1,\cdots,n\},|I|=|J|=d}
 L(c|1^d_u)^{\lambda+\hbar\bar{\epsilon}_J,I}
           _{\lambda\phantom{+\hbar\bar{\epsilon}},I}
 (T_J^\hbar f)(\lambda)
\nonumber
\\
&=&
\sum_{J\subset\{1,\cdots,n\},|J|=d}
\tr_{V(1^d)}
\left(
  \ov   {1^d_{u+c\hbar}}
        {\lambda+\hbar\bar{\epsilon}_J}
        {\lambda}
\circ
  \iv
        {1^d_u}
        {\lambda+\hbar\bar{\epsilon}_J}
        {\lambda}
\right)
 (T_J^\hbar f)(\lambda).
\nonumber
\end{eqnarray}
By the cyclicity property of the trace, the coefficient for $T_J^\hbar$
can be written as
\begin{eqnarray}
\lefteqn{
\tr_{V(1^d)}
\left(
  \ov   {1^d_{u+c\hbar}}
        {\lambda+\hbar\bar{\epsilon}_J}
        {\lambda}
\circ
  \iv
        {1^d_u}
        {\lambda+\hbar\bar{\epsilon}_J}
        {\lambda}
:
V(1^d) \rightarrow
 \path{1^d}
        {\lambda+\hbar\bar{\epsilon}_J}
        {\lambda}
\rightarrow V(1^d)
\right)
}
\nonumber\\
&=&
\tr_{ \path{1^d}
        {\lambda+\hbar\bar{\epsilon}_J}
        {\lambda}
       }
\left(
  \iv
        {1^d_u}
        {\lambda+\hbar\bar{\epsilon}_J}
        {\lambda}
\circ
  \ov   {1^d_{u+c\hbar}}
        {\lambda+\hbar\bar{\epsilon}_J}
        {\lambda}
: \path{1^d}
        {\lambda+\hbar\bar{\epsilon}_J}
        {\lambda}
 \rightarrow
V(1^d) \rightarrow
 \path{1^d}
        {\lambda+\hbar\bar{\epsilon}_J}
        {\lambda}
\right).
\label{eq:trAB=trBA}
\end{eqnarray}
Now taking the trace is much easier because
$\path{1^d}
        {\lambda+\hbar\bar{\epsilon}_J}
        {\lambda}
$
is one dimensional.
We take its basis
(\ref{eq:fusedpath})
and the dual
(\ref{eq:fusedpath*}).
Moreover, we note
\begin{eqnarray*}
  \ov   {1^d_{u+c\hbar}}
        {\lambda+\hbar\bar{\epsilon}_J}
        {\lambda}
(\path{1^d}
        {\lambda+\hbar\bar{\epsilon}_J}
        {\lambda}
)
&=&
  \ov   {\fb_{u+c\hbar}\otimes\cdots\otimes\fb_{u+c\hbar-(d-1)\hbar}}
        {\lambda+\hbar\bar{\epsilon}_J}
        {\lambda}
(\path{1^d}
        {\lambda+\hbar\bar{\epsilon}_J}
        {\lambda}
)
\\
&\subset&
V(\fb_{u+c\hbar}\otimes\cdots\otimes\fb_{u+c\hbar-(d-1)\hbar})
\simeq
V^{\otimes d}
{}.
\end{eqnarray*}
The last isomorphism as vector space is induced by the
obvious identification $V(\fb_u)\simeq V$.
Then we have (Fig.\ref{fig:fcalctr})
\begin{eqnarray}
\lefteqn{
(\ref{eq:trAB=trBA})=
}
\nonumber\\
&&\hspace{-2em}
\tr_{ \path{1^d}
        {\lambda+\hbar\bar{\epsilon}_J}
        {\lambda}
       }
\left(
  \iv
        {1^d_u}
        {\lambda+\hbar\bar{\epsilon}_J}
        {\lambda}
\circ
  \ov   {1^d_{u+c\hbar}}
        {\lambda+\hbar\bar{\epsilon}_J}
        {\lambda}
: \path{1^d}
        {\lambda+\hbar\bar{\epsilon}_J}
        {\lambda}
 \rightarrow
V^{\otimes d}
 \rightarrow
 \path{1^d}
        {\lambda+\hbar\bar{\epsilon}_J}
        {\lambda}
\right)
\nonumber\\
&=&
\tr_{ \path{1^d}
        {\lambda+\hbar\bar{\epsilon}_J}
        {\lambda}
       }
\left(
  \iv
        {\fb_{u}\otimes\cdots\otimes\fb_{u-(d-1)\hbar}}
        {\lambda+\hbar\bar{\epsilon}_J}
        {\lambda}
\circ
  \ov   {\fb_{u+c\hbar}\otimes\cdots\otimes\fb_{u+c\hbar-(d-1)\hbar}}
        {\lambda+\hbar\bar{\epsilon}_J}
        {\lambda}
\right|
_{ \path{1^d}
        {\lambda+\hbar\bar{\epsilon}_J}
        {\lambda}
}
\nonumber\\
&&
\hspace{-2em}
\phantom
{
\tr_{ \path{1^d}
        {\lambda+\hbar\bar{\epsilon}_J}
        {\lambda}
       }
}
\left.
\hspace{5em}
: \path{1^d}
        {\lambda+\hbar\bar{\epsilon}_J}
        {\lambda}
 \rightarrow
V^{\otimes d}
 \rightarrow
 \path{1^d}
        {\lambda+\hbar\bar{\epsilon}_J}
        {\lambda}
\right)
\nonumber\\
&=&
<
{e^*}(1^k)
       ^{\lambda+\hbar\bar{\epsilon}_J}
       _{\lambda}
,
  \iv
        {\fb_{u}\otimes\cdots\otimes\fb_{u-(d-1)\hbar}}
        {\lambda+\hbar\bar{\epsilon}_J}
        {\lambda}
\circ
  \ov   {\fb_{u+c\hbar}\otimes\cdots\otimes\fb_{u+c\hbar-(d-1)\hbar}}
        {\lambda+\hbar\bar{\epsilon}_J}
        {\lambda}
(e(1^k)
       ^{\lambda+\hbar\bar{\epsilon}_J}
       _{\lambda}
 )
>
\nonumber\\
&=&
\sum_{\sigma\in S(d)}
{\rm sgn}(\sigma)
<
e^*_\lambda(j_{1},\cdots,j_{d}),
\nonumber\\
&&
\phantom{
\sum_{\sigma\in S(d)}
}
  \iv
        {\fb_{u}\otimes\cdots\otimes\fb_{u-(d-1)\hbar}}
        {\lambda+\hbar\bar{\epsilon}_J}
        {\lambda}
\circ
  \ov   {\fb_{u+c\hbar}\otimes\cdots\otimes\fb_{u+c\hbar-(d-1)\hbar}}
        {\lambda+\hbar\bar{\epsilon}_J}
        {\lambda}
(
e_\lambda(j_{\sigma(1)},\cdots,j_{\sigma(d)})
)
>
\nonumber\\
&=&
\sum_{\sigma\in S(d)}
{\rm sgn}(\sigma)
\iv
        {{u}}
        {\lambda+\hbar\bar{\epsilon}_{j_1}}
        {\lambda}
\circ
  \ov   {{u+c\hbar}}
        {\lambda+\hbar\bar{\epsilon}_{j_{\sigma(1)}}}
        {\lambda}
\nonumber\\
&&
\phantom{
\sum_{\sigma\in S(d)}
{\rm sgn}(\sigma)
}
\cdots
\iv
        {{u}}
        {\lambda+\hbar\bar{\epsilon}_{j_1}
         +\cdots+\hbar\bar{\epsilon}_{j_{d}}}
        {\lambda+\hbar\bar{\epsilon}_{j_1}
          +\cdots+\hbar\bar{\epsilon}_{j_{d-1}}}
\circ
  \ov   {{u+c\hbar}}
        {\lambda+\hbar\bar{\epsilon}_{j_{\sigma(1)}}
         +\cdots+\hbar\bar{\epsilon}_{j_{\sigma(d)}}}
        {\lambda+\hbar\bar{\epsilon}_{j_{\sigma(1)}}
          +\cdots+\hbar\bar{\epsilon}_{j_{\sigma(d-1)}}}
\nonumber\\
&=&
\sum_{\sigma\in S(d)}
{\rm sgn}(\sigma)
\prod
_{k=1}^d
\left(
\sum_{m=1}^n
\,
\iv
        {{u}}
        {\lambda+\hbar\bar{\epsilon}_{j_1}
         +\cdots+\hbar\bar{\epsilon}_{j_{k}}}
        {\lambda+\hbar\bar{\epsilon}_{j_1}
          +\cdots+\hbar\bar{\epsilon}_{j_{k-1}}}
^{,\,m}
\,
\ov     {{u+c\hbar}}
        {\lambda+\hbar\bar{\epsilon}_{j_{\sigma(1)}}
         +\cdots+\hbar\bar{\epsilon}_{j_{\sigma(k)}}}
        {\lambda+\hbar\bar{\epsilon}_{j_{\sigma(1)}}
          +\cdots+\hbar\bar{\epsilon}_{j_{\sigma(k-1)}}}
_{,\,m}
\right)
{}.
\label{eq:calctrL_d}
\end{eqnarray}
\begin{oekaki}\refstepcounter{oekaki}\addtocounter{oekaki}{-1}
\label{fig:fcalctr}\caption{}\end{oekaki}

\noindent
We are in a position to use the following formula
\begin{eqnarray}
&&\sum_{m=1}^n
\iv
        {u}
        {\mu+\hbar\bar{\epsilon}_j}
        {\mu}
^{,m}
\ov     {u+c\hbar}
        {\lambda+\hbar\bar{\epsilon}_{i}}
        {\lambda}
_{,m}
\nonumber
\\
\nonumber
&=&
\left|
\begin{array}{ccccc}
&(j-1)^{th}&j^{th}&(j+1)^{th}&
\\
\cdots &
\theta_1(\frac{u}{n}-<\mu,\bar{\epsilon}_{j-1}>)
&
\theta_1(\frac{u+c\hbar}{n}-<\lambda,\bar{\epsilon}_{i}>)
&
\theta_1(\frac{u}{n}-<\mu,\bar{\epsilon}_{j+1}>)
& \cdots
\\
&&\vdots&&
\\
\cdots &
\theta_n(\frac{u}{n}-<\mu,\bar{\epsilon}_{j-1}>)
&
\theta_n(\frac{u+c\hbar}{n}-<\lambda,\bar{\epsilon}_{i}>)
&
\theta_n(\frac{u}{n}-<\mu,\bar{\epsilon}_{j+1}>)
& \cdots
\\
&&&&
\end{array}
\right|
\nonumber
\\&&
\quad\cdot
\left|
 \theta_k(\textstyle{\frac{u}{n}}-<\mu,\bar{\epsilon}_{k'}>)
\right|
 ^{-1}
 _{k,k'=1,\cdots,n}
\nonumber
\\
&=&
\frac
    {\theta
     (\frac{c\hbar}{n}+u
     +<\mu,\bar{\epsilon}_j>
     -<\lambda,\bar{\epsilon}_{i}>)}
    {\theta(u)}
\prod_{j'\neq j, 1\leq j'\leq n}
\frac
    {\theta
     (\frac{c\hbar}{n}
     +<\mu,\bar{\epsilon}_{j'}>
     -<\lambda,\bar{\epsilon}_{i}>)}
    {\theta
     (
      <\mu,\bar{\epsilon}_{j'}>
     -<\mu,\bar{\epsilon}_j>)}
,
\label{eq:iv_circ_ov}
\end{eqnarray}
which follows from the definition of $\bar{\phi}$
(\ref{eq:dualityofitv})
and the determinant formula of Vandermonde type
:
\begin{equation}
\mbox{det}
        \left[
         \frac{\theta_{j}(u_k)}
         {\sqrt{-1}\eta(\tau)}
       \right]
        _{j,k=1,\cdots,n}
=
(-1)^{n-1}
\frac
{\theta(\sum_j u_j)}{\sqrt{-1}\eta(\tau)}
\prod_{1\leq j<k \leq n} \frac{\theta(u_k-u_j)}{\sqrt{-1}\eta(\tau)}.
\label{eq:weylkac}
\nonumber
\end{equation}

For convenience, for any statement $P$ let
$$
Y_P :=1 \quad\mbox{\rm if $P$ is true},
\quad
Y_P :=0\quad \mbox{\rm if $P$ is false}.
$$
We have
\begin{eqnarray*}
&&
\sum_{m=1}^n
\iv
        {{u}}
        {\lambda+\hbar\bar{\epsilon}_{j_1}
         +\cdots+\hbar\bar{\epsilon}_{j_{k}}}
        {\lambda+\hbar\bar{\epsilon}_{j_1}
          +\cdots+\hbar\bar{\epsilon}_{j_{k-1}}}
^{,\,m}
\ov     {{u+c\hbar}}
        {\lambda+\hbar\bar{\epsilon}_{j_{\sigma(1)}}
         +\cdots+\hbar\bar{\epsilon}_{j_{\sigma(k)}}}
        {\lambda+\hbar\bar{\epsilon}_{j_{\sigma(1)}}
          +\cdots+\hbar\bar{\epsilon}_{j_{\sigma(k-1)}}}
_{,\,m}.
\\
&=&
\frac
     {\theta
       (
       \frac{c\hbar}{n}+u-(k-1)\hbar
       +<\lambda+\hbar\bar{\epsilon}_{j_1}
          +\cdots+\hbar\bar{\epsilon}_{j_{k-1}}
         ,\bar{\epsilon}_{j_{k}}>
       -<\lambda+\hbar\bar{\epsilon}_{j_{\sigma(1)}}
          +\cdots+\hbar\bar{\epsilon}_{j_{\sigma(k-1)}}
         ,\bar{\epsilon}_{j_{\sigma(k)}}>
     )}
     {\theta(u-(k-1)\hbar)
      }
\\
&&
\cdot
\prod
_{j' \neq j_k, 1\leq j' \leq n}
\frac
     {\theta
       (\frac{c\hbar}{n}
       +<\lambda+\hbar\bar{\epsilon}_{j_1}
          +\cdots+\hbar\bar{\epsilon}_{j_{k-1}}
         ,\bar{\epsilon}_{j'}>
       -<\lambda+\hbar\bar{\epsilon}_{j_{\sigma(1)}}
          +\cdots+\hbar\bar{\epsilon}_{j_{\sigma(k-1)}}
         ,\bar{\epsilon}_{j_{\sigma(k)}}>
     )}
     {\theta
       (
       <\lambda+\hbar\bar{\epsilon}_{j_1}
          +\cdots+\hbar\bar{\epsilon}_{j_{k-1}}
         ,\bar{\epsilon}_{j'}>
       -<\lambda+\hbar\bar{\epsilon}_{j_1}
          +\cdots+\hbar\bar{\epsilon}_{j_{k-1}}
         ,\bar{\epsilon}_{j_{k}}>
     )}
\\
&=&
\frac
       {\theta
       (\frac{c\hbar}{n}+u-(k-1)\hbar
       +<\lambda,\bar{\epsilon}_{j_{k}}
                -\bar{\epsilon}_{j_{\sigma(k)}}
        >
        )}
     {\theta(u-(k-1)\hbar)
      }
\cdot
\prod
_{j' \neq j_k, 1\leq j' \leq n}
\frac
     {\theta
       (\frac{c\hbar}{n}
       +<\lambda,\bar{\epsilon}_{j'}
                -\bar{\epsilon}_{j_{\sigma(k)}}
        >
        +\hbar Y_{j'\in\{j_1,\cdots,j_{k-1}\}}
        )}
     {\theta
       (
       <\lambda,\bar{\epsilon}_{j'}
                -\bar{\epsilon}_{j_{k}}
        >
        +\hbar Y_{j'\in\{j_1,\cdots,j_{k-1}\}}
        )}
\\
&=&
\frac
     {\theta
       (\frac{c\hbar}{n}+u-(k-1)\hbar
       +<\lambda,\bar{\epsilon}_{j_{k}}
                -\bar{\epsilon}_{j_{\sigma(k)}}
        >
        )}
     {\theta(u-(k-1)\hbar)
      }
\\
&&
\cdot
\prod
_{r \neq k, 1\leq r \leq d}
\frac
     {\theta
       (\frac{c\hbar}{n}
       +<\lambda,\bar{\epsilon}_{j_r}
                -\bar{\epsilon}_{j_{\sigma(k)}}
        >
        +\hbar Y_{r<k}
        )}
     {\theta
       (
       <\lambda,\bar{\epsilon}_{j_r}
                -\bar{\epsilon}_{j_k}
        >
        +\hbar Y_{r<k}
        )}
\prod_{j'\notin J}
\frac
     {\theta
       (\frac{c\hbar}{n}
       +<\lambda,\bar{\epsilon}_{j'}
                -\bar{\epsilon}_{j_{\sigma(k)}}
        >
        )}
     {\theta
       (
       <\lambda,\bar{\epsilon}_{j'}
                -\bar{\epsilon}_{j_k}
        >
        )}
{}.
\end{eqnarray*}
Therefore
\begin{eqnarray}
\lefteqn{
(\ref{eq:calctrL_d})
}
\nonumber
\\
&=&
\sum
_{\sigma\in S(d)}
{\rm sgn}(\sigma)
\prod
_{k=1}^d
\left(
\frac
     {\theta
       (\frac{c\hbar}{n}+u-(k-1)\hbar
       +<\lambda,\bar{\epsilon}_{j_{k}}
                -\bar{\epsilon}_{j_{\sigma(k)}}
        >
        )}
     {\theta(u-(k-1)\hbar)
      }
\right.
\nonumber
\\
&&
\left.
\phantom
 {
 \sum
 _{\sigma\in S(d)}
 {\rm sgn}(\sigma)
 \prod
 _{k=1}^d
 }
\cdot
\prod
_{r \neq k, 1\leq r \leq d}
\frac
     {\theta
       (\frac{c\hbar}{n}
       +<\lambda,\bar{\epsilon}_{j_r}
                -\bar{\epsilon}_{j_{\sigma(k)}}
        >
        +\hbar Y_{r<k}
        )}
     {\theta
       (
       <\lambda,\bar{\epsilon}_{j_r}
                -\bar{\epsilon}_{j_k}
        >
        +\hbar Y_{r<k}
        )}
\prod_{j'\notin J}
\frac
     {\theta
       (\frac{c\hbar}{n}
       +<\lambda,\bar{\epsilon}_{j'}
                -\bar{\epsilon}_{j_{\sigma(k)}}
        >
        )}
     {\theta
       (
       <\lambda,\bar{\epsilon}_{j'}
                -\bar{\epsilon}_{j_k}
        >
        )}
\right)
\nonumber
\\
&=&
\left(
\frac
     {
\sum
_{\sigma\in S(d)}
{\rm sgn}(\sigma)
\prod
_{k=1}^d
\prod
_{r=1}^d
\theta
       (\frac{c\hbar}{n}
       +<\lambda,\bar{\epsilon}_{j_r}
                -\bar{\epsilon}_{j_{\sigma(k)}}
        >
        +\hbar Y_{r<k}
        +\delta_{k,r}(u-(k-1)\hbar)
        )
}
{
\prod_{r=0}^{d-1}\theta(u-r\hbar)
\prod_{1\leq r<k \leq d}
       \theta
       (
       <\lambda,\bar{\epsilon}_{j_r}
                -\bar{\epsilon}_{j_k}
        >
        +\hbar
        )
\prod_{1\leq k<r \leq d}
       \theta
       (
       <\lambda,\bar{\epsilon}_{j_r}
                -\bar{\epsilon}_{j_k}
        >
        )
}
\right)
\nonumber\\
&&
\cdot
\prod_{j\in J, j'\notin J}
\frac
     {\theta
       (\frac{c\hbar}{n}
       +<\lambda,\bar{\epsilon}_{j'}
                -\bar{\epsilon}_{j}
        >
        )}
     {\theta
       (
       <\lambda,\bar{\epsilon}_{j'}
                -\bar{\epsilon}_{j}
        >
        )}
\nonumber\\
&=&
\frac
     {{\rm det}\left[
\prod
_{r=1}^d
\theta
       (\frac{c\hbar}{n}
       +<\lambda,\bar{\epsilon}_{j_r}
                -\bar{\epsilon}_{j_{k'}}
        >
        +\hbar Y_{r<k}
        +\delta_{k,r}(u-(k-1)\hbar)
        )
\right]_{k,k'=1,\cdots,d}}
{
\prod_{r=0}^{d-1}\theta(u-r\hbar)
\prod_{1\leq r<k \leq d}
       \theta
       (
       <\lambda,\bar{\epsilon}_{j_r}
                -\bar{\epsilon}_{j_k}
        >
        +\hbar
        )
\prod_{1\leq k<r \leq d}
       \theta
       (
       <\lambda,\bar{\epsilon}_{j_r}
                -\bar{\epsilon}_{j_k}
        >
        )
}
\nonumber\\
&&
\cdot
\prod_{j\in J, j'\notin J}
\frac
     {\theta
       (\frac{c\hbar}{n}
       +<\lambda,\bar{\epsilon}_{j'}
                -\bar{\epsilon}_{j}
        >
        )}
     {\theta
       (
       <\lambda,\bar{\epsilon}_{j'}
                -\bar{\epsilon}_{j}
        >
        )}.
\label{eq:calctrLend}
\end{eqnarray}
\noindent
To end the proof,
we are led to show the following Lemma \ref{lem:qFay} below,
which is quite interesting itself.
According to the lemma, we see the unwanted factor in
(\ref{eq:calctrLend}) cancels as desired. QED.
\hfill\fb

\medskip

\begin{lem}
Recall $Y_{r<s} :=1$ if $r<s$ holds, $Y_{r<s}:=0$ otherwise.
The following formula holds:
\nobreak
\begin{eqnarray}
&&
{\rm det}\left[
 \prod_{r=1}^d
  \theta
    \left(\mu_r - \lambda_{s'} +\hbar Y_{r<s} + \delta_{r,s}(u-(s-1)\hbar)
    \right)
\right]_{s,s'=1,\cdots,d}
\label{eq:qFay}
\\
\nopagebreak
&=&
\nonumber
\theta(u+\sum_{r=1}^d(\mu_r-\lambda_r))
\prod_{s=1}^{d-1}\theta(u-s\hbar)
\prod_{1\leq s<s'\leq d}
 \theta(\lambda_{s'}-\lambda_{s})
 \theta(\hbar + \mu_s-\mu_{s'}).
\end{eqnarray}
\label{lem:qFay}
\end{lem}
\noindent
Here the variables $\{\lambda_s\}, \{\mu_s\}$ are arbitrary, i.e.
neither the conditions
$\sum_{s=1}^d \lambda_s =0$ nor $\sum_{s=1}^d \mu_s=0$
are necessary.
Proof of this formula uses the induction on $d$
and will be given in the Appendix.
The $\hbar = 0$ case of (\ref{eq:qFay}) is easily transformed into
the Cauchy type determinant formula
\begin{equation}
{\rm det}\left[
\frac{\theta(\mu_s-\lambda_{s'}+u)}
     {\theta(\mu_s-\lambda_{s'})\theta(u)}
\right]_{s,s'=1,\cdots,d}
=
\frac
{\theta(u+\sum_{r=1}^d(\mu_r-\lambda_r))}
{\theta(u)}
\frac
{
\prod_{1\leq s<s'\leq d}
 \theta(\mu_s-\mu_{s'})
 \theta(\lambda_{s'}-\lambda_{s})
}
{
\prod_{s,s'=1,\cdots,d}
\theta(\mu_s-\lambda_{s'})
},
\label{eq:Fay}
\end{equation}
which is also known as
the genus 1 case of Fay's trisecant formula \cite{Fay}.
But I do not know whether
(\ref{eq:qFay}) is previously known or not.

It is also interesting to remark that this $\hbar =0$ case was
quite relevant in \cite{Ruij} although his approach for the commuting
system is different from ours.


\section{Relations to other approaches for the system}

\subsection{Ruijsenaars' operators}
\begin{prop}
Put $g:=\frac{c}{n}$.
Define a function $\Phi$ on ${\bf h}^*$ by
\begin{equation}
\Phi(\lambda) :=
\prod_{k\neq k'} d^+(z_k/z_{k'}),
\quad
d^+(z):=
\prod_{k=0}^{\infty} \prod_{m=0}^\infty
 \frac {1-zq^{m+1}p^k} {1-zq^{m+g+1}p^k}
 \frac {1-z^{-1}q^{m-g}p^{k+1}} {1-z^{-1}q^m p^{k+1}},
\label{eq:d^+}
\end{equation}
where
$p=\exp 2\pi i \tau, q=\exp 2\pi i \hbar (|q|<1)$
and
$z_j := \exp 2\pi i <\lambda, \bar{\epsilon}_j>$.
Then the conjugation by the square root $\Phi^{1/2}$ yields
\footnote {as long as it makes sense}
Ruijsenaars' \mbox{\rm\cite{Ruij}} commuting operators:
\begin{eqnarray*}
&&
\left(\frac{\theta(v+\frac{dc\hbar}{n})}{\theta(v)}\right)^{-1}
\cdot
\Phi^{-1/2} M_{d}(c|u) \Phi^{1/2}
\\
&=&
\sum_
     {I\subset\{1,\cdots,n\},|I|=d}
 \left(\prod
      _{s \notin I, t \in I}
  \sqrt{
        \frac{\theta(\frac{c}{n}\hbar+<\lambda, \epsilon_s-\epsilon_t>)}
             {\theta(<\lambda, \epsilon_s-\epsilon_t>)}
        }\right)
T_I^\hbar
 \left(\prod
      _{s \notin I, t \in I}
  \sqrt{
        \frac{\theta(\frac{c}{n}\hbar+<\lambda, \epsilon_t-\epsilon_s>)}
             {\theta(<\lambda, \epsilon_t-\epsilon_s>)}
        }\right).
\end{eqnarray*}
\label{prop:reltoruij}
\end{prop}

{\bf Proof.}[cf.\cite{D}]
The function $d^+$ (\ref{eq:d^+}) satisfies

$$
\frac{d^+(z)}{d^+(qz)}
=
\prod_{k=0}^{\infty}
 \frac {1-zq p^k}{1-z q^{1+g}p^{k}}
 \frac {1-(zq)^{-1}p^{k+1}} {1-(zq^{1+g})^{-1}p^{k+1}}
$$
and consequently $\Phi$ satisfies
\begin{eqnarray*}
&&\frac{\Phi(z)}{T^\hbar_i\Phi(z)}
=
\frac{\Phi(z_1,\cdots,z_n)}{\Phi(z_1,\cdots,qz_i,\cdots,z_n)}
\\
&=&
\prod_{k=0}^\infty
\prod_{j\neq i}
\frac
 {(1-\frac{qz_i}{z_j}p^k)
   (1-\frac{z_j}{qz_i}p^{k+1})
    (1-\frac{z_j}{z_i}q^gp^k)
     (1-\frac{z_i}{z_j}q^{-g}p^{k+1})}
 {(1-\frac{qz_i}{z_j}q^gp^k)
   (1-\frac{z_j}{qz_i}q^{-g}p^{k+1})
    (1-\frac{z_j}{z_i}p^k)
     (1-\frac{z_i}{z_j}p^{k+1})}
\\
&=&
\prod_{j\neq i}\left\{
\frac
{(1-\frac{qz_i}{z_j})(1-\frac{z_j}{z_i}q^g)}
{(1-\frac{qz_i}{z_j}q^g)(1-\frac{z_j}{z_i})}
\prod_{k=1}^\infty
\frac
 {(1-(\frac{qz_i}{z_j})p^k)
   (1-(\frac{qz_i}{z_j})^{-1}p^k)
    (1-(\frac{z_j}{z_i}q^g)p^k)
     (1-(\frac{z_j}{z_i}q^g)^{-1}p^k)}
 {(1-(\frac{qz_i}{z_j}q^g)p^k)
   (1-(\frac{qz_i}{z_j}q^g)^{-1}p^k)
    (1-\frac{z_j}{z_i}p^k)
     (1-(\frac{z_j}{z_i})^{-1}p^k)}
\right\}
\\
&=&
\prod_{j\neq i}
\frac
{\theta(\hbar+\lambda_{i,j})
  \theta(g\hbar+\lambda_{j,i})}
{\theta(g\hbar+\hbar+\lambda_{i,j})
  \theta(\lambda_{j,i})},
\end{eqnarray*}
where $\lambda_{i,j}=<\lambda,\bar{\epsilon}_i-\bar{\epsilon}_j>$.
 From this property, for general $I\subset\{1,\cdots,n\}$
we have
\begin{eqnarray*}
\frac{\Phi(z)}{T^\hbar_I\Phi(z)}
&=&
\prod_{i\in I,j\notin I}
\frac
{\theta(\hbar+\lambda_{i,j})
  \theta(g\hbar+\lambda_{j,i})}
{\theta(g\hbar+\hbar+\lambda_{i,j})
  \theta(\lambda_{j,i})}.
\end{eqnarray*}
Therefore,
\begin{eqnarray*}
\lefteqn{
\left(\frac{\theta(v+\frac{dc\hbar}{n})}{\theta(v)}\right)^{-1}
\cdot
\left(
\Phi^{-1/2} M_d(c|u) \Phi^{1/2} f\right)(\lambda)
}
\\
&=&
\sum_{|I|=d}
\left(
 \prod_{i\in I,j\notin I}
  \frac{\theta(\lambda_{j,i}+g\hbar)}
      {\theta(\lambda_{j,i})}
\right)
\cdot
\left(
\Phi^{-1/2}(T_I^\hbar \Phi)^{1/2}
\right)(\lambda)
(T_I^\hbar f)(\lambda)
\\
&=&
\sum_{|I|=d}
\left(
 \prod_{i\in I,j\notin I}
  \frac{\theta(\lambda_{j,i}+g\hbar)}
      {\theta(\lambda_{j,i})}
\right)
\left(
\prod_{i\in I,j\notin I}
\frac
{\theta(g\hbar+\hbar+\lambda_{i,j})
  \theta(\lambda_{j,i})}
{\theta(\hbar+\lambda_{i,j})
  \theta(g\hbar+\lambda_{j,i})}
\right)^{1/2}
(T_I^\hbar f)(\lambda)
\\
&=&
\sum_{|I|=d}
\left(
 \prod_{i\in I,j\notin I}
  \frac{\theta(\lambda_{j,i}+g\hbar)}
      {\theta(\lambda_{j,i})}
\right)^{1/2}
\left(
\prod_{i\in I,j\notin I}
\frac
{\theta(g\hbar+\hbar+\lambda_{i,j})}
{\theta(\hbar+\lambda_{i,j})}
\right)^{1/2}
(T_I^\hbar f)(\lambda)
\\
&=&
\sum_{|I|=d}
\left(
 \prod_{i\in I,j\notin I}
  \frac{\theta(\lambda_{j,i}+g\hbar)}
      {\theta(\lambda_{j,i})}
\right)^{1/2}
T_I^\hbar
\left(
\prod_{i\in I,j\notin I}
\frac
{\theta(g\hbar+\lambda_{i,j})}
{\theta(\lambda_{i,j})}
\right)^{1/2}
f.
\end{eqnarray*}
\hfill \fb

{\bf Remark}.
Put
$
d(z):=d_+(z/q), \delta :=\prod_{k\neq k'}d(z_k/z_{k'})
$
and define $\bar{f}(\lambda):=f(-\lambda)$.
Consider the inner product defined by
$
(f,g)':= [\bar{f}g\delta]_1,
$
where $[\phantom{fg}]_1$ stands for the constant term in $\lambda$.
Then generalizing the trigonometric (Macdonald's) case,
one can directly verify that the operators $\{M_d(c|u)\}$ are
formally self-adjoint with respect to this inner product.

\subsection{Krichever's Lax matrix}

Put
$
\dot{M}_d:=\frac{\theta(u)}{\theta(u+\frac{cd}{n}\hbar)}M_d.
$
In the limit $\hbar \rightarrow 0$,
one can check that our system degenerates to the
elliptic Calogero-Moser system
(for the proof see (\ref{eq:l'Hospital}) in
Subsection \ref{sect:difflim}),
\begin{equation}
\frac{1}{\hbar^2}
\left(
-2\dot{M}_2+\dot{M}_1^2-2\dot{M}_1+n
\right)
\stackrel{\hbar\rightarrow 0}{\longrightarrow}
H,
\label{eq:limisOP}
\end{equation}
$$
H = \Delta^{c/n}
\circ
\left(
\sum_{i=1}^n
 \frac{\partial^2}
 {\partial\bar{\epsilon}_i^2}
-
\frac{c}{n}(\frac{c}{n}+1)
\sum_{i<j}
   (\log{\theta})''(\lambda_{i,j})
\right)
\circ
\Delta^{-c/n}.
$$
Here
$
\Delta(\lambda):=\prod_{i<j}\theta(\lambda_{i,j}).
$
Note that the Jacobi $\theta$ function and the Weierstrass sigma
function with quasi-period $1,\tau$ are in the relation
$
\sigma(u) = e^{hu^2}{\theta(u)}/{\theta'(0)}
$,
$
h =
-2\pi\sqrt{-1}
\frac{\partial}{\partial\tau}\log{\eta(\tau)}
$
so that we have
$
\wp(u)
 = -(\log{\sigma(u)})''
 = -(\log{\theta(u)})''+ h
$
\cite{CH}.
It is well known that this system is integrable
and admits a Lax matrix
formalism \cite{Krichever}.
Here we shall derive the Lax matrix from our L-operator.

It turns out that the L-operator (\ref{eq:defofL})
should be slightly modified
for our purpose.
We conjugate it by the matrix formed by the intertwining
vectors:
$$
L(c|u)^i_j \rightarrow
\tilde{L}(c|u)^{i'}_{j'}
:=
\sum_{i,j=1}^n
\iv     {u}
        {\lambda+\hbar\bar{\epsilon}_{j'}}
        {\lambda}
        ^{j}
\ov     {u}
        {\lambda+\hbar\bar{\epsilon}_{i'}}
        {\lambda}
        _{i}
L(c|u)^i_j
{}.
$$
Using the property (\ref{eq:dualityofitv}), this conjugated matrix
takes the simple form (Fig. \ref{fig:fLtilde})
$$
\tilde{L}(c|u)^{i}_{j}
=
\sum_{k=1}^n
\iv     {u}
        {\lambda+\hbar\bar{\epsilon}_{j}}
        {\lambda}
        ^{k}
\ov     {u+c\hbar}
        {\lambda+\hbar\bar{\epsilon}_{i}}
        {\lambda}
        _{k}
T_i^\hbar.
$$
\begin{oekaki}\refstepcounter{oekaki}\addtocounter{oekaki}{-1}
\label{fig:fLtilde}\caption{}\end{oekaki}
\noindent
Moreover, for the summation in the right hand side we can apply the
formula (\ref{eq:iv_circ_ov}):
\begin{eqnarray*}
\sum_{k=1}^n
\iv     {u}
        {\lambda+\hbar\bar{\epsilon}_{j}}
        {\lambda}
        ^{k}
\ov     {u+c\hbar}
        {\lambda+\hbar\bar{\epsilon}_{i}}
        {\lambda}
        _{k}
&=&
\frac{\theta(\frac{c\hbar}{n}+u+\lambda_{j,i})}{\theta(u)}
\prod_{k\neq j}
\frac{\theta(\frac{c\hbar}{n}+\lambda_{k,i})}{\theta(\lambda_{k,j})},
\end{eqnarray*}
so that we have the formula
\begin{equation}
\tilde{L}(c|u)_j^i
=
\left(
\frac{\theta(\frac{c\hbar}{n}+u+\lambda_{j,i})}{\theta(u)}
\prod_{k\neq j}
\frac{\theta(\frac{c\hbar}{n}+\lambda_{k,i})}{\theta(\lambda_{k,j})}
\right)
T_i^\hbar.
\label{eq:tildeLformula}
\end{equation}
\noindent
This is our version of the Lax matrix suggested in \cite{Ruij}.
 From this it is easy to see:
\begin{equation}
\tilde{L}(c|u)_j^i
\stackrel{\hbar\rightarrow 0}{\longrightarrow}
\delta_j^i.
\label{eq:limL=id}
\end{equation}
To recover Krichever's Lax matrix, we look at the
coefficient for $\hbar^1$ in $\tilde{L}$.
We have
\begin{eqnarray*}
&&
\left.\frac{\partial}{\partial\hbar}
\tilde{L}(c|u)^i_j
\right|_{\hbar=0}
=
\left.\frac{\partial}{\partial\hbar}
\left(
\frac{\theta(\frac{c\hbar}{n}+u+\lambda_{j,i})}{\theta(u)}
\prod_{k\neq j}
\frac{\theta(\frac{c\hbar}{n}+\lambda_{k,i})}{\theta(\lambda_{k,j})}
\right)
T_i^\hbar
\right|_{\hbar=0}
\\
&=&
\left.\frac{\partial}{\partial\hbar}
\left(
\frac{\theta(\frac{c\hbar}{n}+u+\lambda_{j,i})}{\theta(u)}
\prod_{k\neq j}
\frac{\theta(\frac{c\hbar}{n}+\lambda_{k,i})}{\theta(\lambda_{k,j})}
\right)
\right|_{\hbar=0}
\cdot id
+
\delta_j^i \frac{\partial}{\partial\bar{\epsilon}_i}
\\
&=&
\frac{c}{n}
\left\{
\frac{\theta'(u+\lambda_{j,i})}{\theta(u)}\delta_j^i
+
\frac{\theta(u+\lambda_{j,i})}{\theta(u)}
\left(
 \delta_j^i\sum_{k\neq j}
  \frac{\theta'(\lambda_{k,j})}{\theta(\lambda_{k,j})}
 +
 (1-\delta_j^i)
  \theta'(0)
  \frac
   {\prod_{k\neq i,j}\theta(\lambda_{k,i})}
   {\prod_{k\neq j}\theta(\lambda_{k,j})}
\right)
\right\}
+
\delta_j^i \frac{\partial}{\partial\bar{\epsilon}_i}
\\
&=&
\delta_j^i
\left\{\frac{c}{n}
 \left(
  \frac{\theta'(u)}{\theta(u)}
  - 
  \frac{\partial}{\partial\bar{\epsilon}_i}
   (\log \Delta(\lambda))
\right)
+
\frac{\partial}{\partial\bar{\epsilon}_j}
\right\}
+
 (1-\delta_j^i)
  \frac
   {\theta(u+\lambda_{j,i})\theta'(0)}
   {\theta(u)\theta(\lambda_{j,i})}
  \frac
   {\prod_{k\neq i}\theta(\lambda_{k,i})}
   {\prod_{k\neq j}\theta(\lambda_{k,j})}
{}.
\end{eqnarray*}
If we further modify this matrix by
1) the conjugation by $\Delta(\lambda)^{c/n}$
and by
2) the similarity transformation by the diagonal matrix
$
diag
\left(
{\prod_{k\neq i}\theta(\lambda_{k,i})}
\right)_{i=1,\cdots,n},
$
we obtain
\begin{prop}
\begin{eqnarray*}
&&
  \frac
   {\prod_{k\neq j}\theta(\lambda_{k,j})}
   {\prod_{k\neq i}\theta(\lambda_{k,i})}
\cdot
\Delta^{-c/n}
\circ
\left.\frac{\partial}{\partial\hbar}
\tilde{L}(c|u)^i_j
\right|_{\hbar=0}
\circ
\Delta^{c/n}
\\
&&
=
\delta_j^i
\left\{\frac{c}{n}
  \frac{\theta'(u)}{\theta(u)}
+
\frac{\partial}{\partial\bar{\epsilon}_j}
\right\}
+
 (1-\delta_j^i)
\frac{c}{n}
  \frac
   {\theta(u+\lambda_{j,i})\theta'(0)}
   {\theta(u)\theta(\lambda_{j,i})}
=: K(c|u)_j^i          
{}.
\end{eqnarray*}
\label{prop:KricheverL}
\end{prop}
One sees that this matrix $K(c|u)$
essentially coincides with Krichever's Lax matrix modulo
the scalar matrix addition
$\frac{c}{n}\frac{\theta'(u)}{\theta(u)}id.$

\medskip
In the classical case,
the so-called ``dynamical r-structure'' for Krichever's Lax matrix
is obtained by Sklyanin \cite{Skldr}.
 From our viewpoint, this can be easily extended to the quantum case
and the dynamical r-matrix can be essentially
understood as the first derivative of the face weight at $\hbar=0$.
Thus the dynamical variable dependence of the dynamical r-matrix
quite naturally arise since the face weight depend on the variable
 $\lambda \in {\bf h}^*$.
We will treat these points in a seperate paper.

It should be also mentioned that the matrix $\tilde{L}(c|u)$ gives
an representation of the quadratic algebra in \cite{Felder} after a
suitable formulation.
It is known \cite{BD} that the
elliptic extension of the vertex type solution to the Yang-Baxter
equation only exist for the type A case (i.e. the Belavin R-matrix),
while the face type solutions (face weights)
are known to exist for the other types as well \cite{JMO2}.
This may suggest that, to generalize the Olshanetsky-Perelomov system
\cite{OP}
with the symmetry other than type A,
searching the analogue of the $\tilde{L}(c|u)$ matrix might be useful
(cf.\cite{Avanetal}).

\subsection{The generating function and the differential limit}

Conceptually, our commuting operators $\{M_d(c|u)\}_{d=1}^n$ are
generated by the L-operator $L(c|u)=[L(c|u)_j^i]$ by fusing and taking
the trace. In this sence $L(c|u)$ is the generating function of the
system. 
More explicitly,
we can state the relation in the following familiar form:

\begin{thm}
$$
\sum_{d=0}^n (-t)^{n-d} M_d(c|u) = {\np}{\rm det}[L(c|u)-t]{\np},
$$
where we put $M_0(c|u):={\rm id}$
and the normal product {\np\np} is defined by
putting the difference/differential operators to the right:
e.g.
$
\np f(\lambda) T^\hbar_I g(\lambda)T^\hbar_J \np
:= f(\lambda)g(\lambda) T^\hbar_I T^\hbar_J
$
for any function $f, g$ in $\lambda$.
\label{thm:genfnc}
\end{thm}

{\bf Proof.}
As in the trigonometric (Macdonald's) case,
we can rewrite our operator as follows:
$$
{M}_d(c|u)
=
\sum_{|I|=d}
\frac{\theta(u+\frac{cd}{n}\hbar)}{\theta(u)}
\frac{T^{-\frac{c}{n}\hbar}_I\Delta}{\Delta}(\lambda)
T^\hbar_I,
$$
where
$\Delta(\lambda)=\prod_{i<j}\theta(\lambda_j-\lambda_i).$
In view of the determinant formula (\ref{eq:weylkac}),
we have
$$
\frac{\theta(u+\frac{cd}{n}\hbar)}{\theta(u)}
\frac{T^{-\frac{c}{n}\hbar}_I\Delta}{\Delta}(\lambda)
=
\frac{      T^{-\frac{c}{n}\hbar}_I
      {\rm det}[\theta_j(\frac{u}{n}-\lambda_i)]}
{{\rm det}[\theta_j(\frac{u}{n}-\lambda_i)]}
\;\;
\mbox{for $I\subset\{1,\cdots,n\},|I|=d$}.
$$
Therefore,
\begin
{eqnarray}
\sum_{d=0}^n (-t)^{n-d} M_d(c|u)
&=&
\sum_{d=0}^n (-t)^{n-d}
\sum_{|I|=d}
\left(
\frac{      T^{-\frac{c}{n}\hbar}_I
      {\rm det}[\theta_j(\frac{u}{n}-\lambda_i)]}
{{\rm det}[\theta_j(\frac{u}{n}-\lambda_i)]}
\right)
T^\hbar_I
\nonumber\\
&=&
\sum_I (-t)^{n-|I|}
\frac
{{\np}{\rm det}
 [(T^{-\frac{c}{n}\hbar Y_{i\in I}}_i
  \theta_j(\frac{u}{n}-\lambda_i)
  )
   T^{\hbar Y_{i \in I}}_i]{\np}}
{{\rm det}[\theta_j(\frac{u}{n}-\lambda_i)]}
\nonumber\\
&=&
\frac{
{\np}{\rm det}[
\theta_j(\frac{u+c\hbar}{n}-\lambda_i) T^{\hbar}_i
-t\cdot \theta_j(\frac{u}{n}-\lambda_i)]{\np}}
{{\rm det}[\theta_j(\frac{u}{n}-\lambda_i)]}
\label{eq:Sekiguchi}\\
&=&
{\np}
{\rm det}\left[
\theta_j(\textstyle\frac{u+c\hbar}{n}-\lambda_i) T^{\hbar}_i
-t\cdot \theta_j(\textstyle\frac{u}{n}-\lambda_i)\right]
\cdot
{{\rm det}\left[\theta_j(\textstyle\frac{u}{n}-\lambda_i)\right]}^{-1}
{\np}
\label{eq:multfromR}
\\
&=&
{\np}{\rm det}
\left[
\sum_i
\phi(u+c\hbar)_{\lambda,j}^{\lambda+\hbar\bar{\epsilon_i}}
 T^{\hbar}_i
\bar{\phi}(u)_{\lambda}^{\lambda+\hbar\bar{\epsilon_i},k}
-t\delta_j^k
\right]_{j,k=1,\cdots,n}
{\np}
\nonumber
\\
&=&
{\np}{\rm det}[L(c|u)-t\,]{\np}.
\nonumber
\label{eq:MviaL}
\end{eqnarray}
\hfill\fb

{\bf Remark.}
Let
${\bf 1}:=\sum_{i=1}^n \epsilon_i$.
Then
$[M_d(c|u),T_{\bf 1}^{a\hbar}]=0$
for any $a$.
This ensures that the replacement of
$T_i^{a\hbar}= \exp{a\hbar\frac{\partial}{\partial\bar{\epsilon}_i}}$
by
$\exp{a\hbar\frac{\partial}{\partial\epsilon_i}}$ in the formula of
$M_d(c|u)$
will not affect the commuting nature of the family.
Thus we can and will abuse the symbol
$T_i^{a\hbar}$ (resp. $\partial_i$)
for these two meanings,
$\exp{a\hbar\bar{\epsilon}_i}$ and $\exp{a\hbar{\epsilon}_i}$
(resp.
$\frac{\partial}{\partial\bar{\epsilon}_i}$
and
$\frac{\partial}{\partial{\epsilon}_i}$),
in what follows.
\hfill \fb

\medskip
The numerator in (\ref{eq:Sekiguchi}) can be regarded as the elliptic
difference analogue of Jiro Sekiguchi's generating function for
the trigonometric differential case \cite{JiroSekiguchi}.

In (\ref{eq:multfromR}), there is another option, i.e.
multipling
 ${{\rm det}\left[\theta_j(\textstyle\frac{u}{n}-\lambda_i)\right]}^{-1}$
from the left.
As is easily seen, the result is
\begin{equation}
\sum_{d=0}^n (-t)^{n-d} M_d(c|u)
=
{\np}{\rm det}
\left[
\sum_j
\bar{\phi}(u)_{\lambda}^{\lambda+\hbar\bar{\epsilon_k},j}
\phi(u+c\hbar)_{\lambda,j}^{\lambda+\hbar\bar{\epsilon_i}}
 T^{\hbar}_i
-t\delta_i^k
\right]_{i,k=1,\cdots,n}{\np}
=
{\np}{\rm det}[\tilde{L}(c|u)-t\,]{\np}
,
\label{eq:MviaLtilde}
\end{equation}
where $\tilde{L}$ is given by (\ref{eq:tildeLformula}).
Using this generating function, we can derive the commuting operators
in the differential limit $\hbar\rightarrow 0$.
Since
$
\tilde{L}
\stackrel{\hbar\rightarrow 0}{\rightarrow}
id
$
(\ref{eq:limL=id}),
$
\frac{1}{\hbar}(\tilde{L}-1)
\stackrel{\hbar\rightarrow 0}{\rightarrow}
\left.\frac{\partial}{\partial\hbar}\tilde{L}_\hbar(c|u)\right|_{\hbar=0}
=:\tilde{L}'_0(c|u),
$
we have 
$$
\np{\rm det}[\tilde{L}_\hbar(c|u)-(1+\hbar t)\,]\np
=
\np{\rm det}[\hbar(\tilde{L}'_0(c|u)-t)+o(\hbar)]\np
\stackrel{\hbar\rightarrow 0}{\sim}
\hbar^n\np{\rm det}[\tilde{L}'_0(c|u)-t]\np+o(\hbar^n).
$$
This leads us to define the differential operators
$\{D_d(c|u)\}$ by
\begin{equation}
\lim_{\hbar\rightarrow 0}
\np
{\rm det}\left[\frac{\tilde{L}_\hbar(c|u)-(1+\hbar t)}{\hbar}\right]
\np
=
\np{\rm det}[\tilde{L}'_0(c|u)-t]\np
=: \sum_{d=0}^n(-t)^{n-d}D_d(c|u).
\label{eq:defD_d}
\end{equation}
In terms of Krichever's Lax matrix (Prop. 3), we have
an equivalent definition
$$
\Delta^{c/n}\circ\np{\rm det}[K(c|u)-t]\np\circ\Delta^{-c/n}
= \sum_{d=0}^n(-t)^{n-d}D_d(c|u).
$$
These formula actually define a commuting family
$[D_d(c|u),D_{d'}(c|v)]=0$,
because the left hand sides are just the summation of the commuting
operators $\{M_d\}$.
Obviously, $D_0=id$.
The higher operators are given by:
\begin{eqnarray*}
\np{\rm det}[\tilde{L}'_0-t]\np
&=&
\sum_{d=0}^n (-t)^{n-d}
\sum_{I\subset \{1,\cdots,n\},|I|=d}
\np{\rm det}[\tilde{L}'_0(c|u)^i_{i'}]_{i,i'\in I}\np
\\
&=&
\sum_{d=0}^n (-t)^{n-d}
\sum_{|I|=d}
\sum_{J\subset I}
\;
{\rm det}[C'_0(u,\lambda)^j_{j'}]_{j,j'\in J}
\cdot
\partial^{I\backslash J},
\end{eqnarray*}
where
$
\partial^I :=\prod_{i\in I}\partial_i
$
and
$$
C'_0(u,\lambda)^i_{j}
:=
\left.
\frac{\partial}{\partial\hbar}{C_\hbar(u,\lambda,\lambda)^j_{j'}}
\right|_{\hbar=0},
\;
C_\hbar(u,\lambda,\mu)^i_{j}
:=
\frac{\theta(\frac{c\hbar}{n}+u+\mu_j-\lambda_{i})}{\theta(u)}
\prod_{k\neq j}
\frac{\theta(\frac{c\hbar}{n}+\mu_k-\lambda_i)}{\theta(\mu_{k,j})}
\;\;
\mbox{(cf.(\ref{eq:tildeLformula}))}.
$$
Since
$
C'_0(u,\lambda)^i_{j}
=
-\textstyle\frac{c}{n}
\left.
\frac{\partial}{\partial{\lambda_i}}C_0(u,\lambda,\mu)^i_{j}
\right|_{\lambda=\mu}
,
$
we can apply the formula (\ref{eq:Fay}) as follows:
\begin{eqnarray}
{\rm det}[C'_0(u)^j_{j'}]_{j,j'\in J}
&=&
(-\textstyle\frac{c}{n})^{|J|}
{\rm det}
\left[\left.
\frac{\partial}{\partial{\lambda_i}}C_0(u,\lambda,\mu)
\right|_{\lambda=\mu}
\right]_{i,j\in J}
\nonumber
\\
&=&
(-\textstyle\frac{c}{n})^{|J|}
\left.
(
\prod_{i\in J}\frac{\partial}{\partial{\lambda_i}}
)
{\rm det}
\left[
C_0(u,\lambda,\mu)
\right]_{i,j\in J}
\right|_{\lambda=\mu}
\nonumber
\\
&\stackrel{(\ref{eq:Fay})}{=}&
(-\textstyle\frac{c}{n})^{|J|}
\left(
\prod_{i\in J}\frac{\partial}{\partial{\lambda_i}}
\right)
\left.
\frac
{\theta(u+\sum_{j\in J}(\mu_j-\lambda_j))\Delta_J(\lambda)\Delta_J(-\mu)}
{\theta(u)\Delta_J(\mu)\Delta_J(-\mu)}
\right|_{\lambda=\mu}
\\
&=&
(-\textstyle\frac{c}{n})^{|J|}
\displaystyle\sum_{J=J'\sqcup J''}(-1)^{|J'|}
\frac{\theta^{(|J'|)}(u)}{\theta(u)}
\frac{\partial^{J''}\Delta_J}{\Delta_J}(\lambda).
\nonumber
\end{eqnarray}
Here
$\Delta_J(\lambda)
:=\prod_{j,j'\in J, j<j'}\theta(\lambda_j-\lambda_{j'})$
{}.
Summarizing, 
(note that
$\frac{\partial^{J''}\Delta_J}{\Delta_J}
=\frac{\partial^{J''}\Delta}{\Delta}$
for
$J''\subset J \subset \{1,\cdots,n\}$)
\begin{prop}
Put $(-\frac{n}{c}\partial)^J:=\prod_{j\in J}-\frac{n}{c}\partial_j$.
In the differential limit $\hbar\rightarrow 0$, the commuting
operators
$\{D_d(c|u)\}_{d=1}^{n}$ (\ref{eq:defD_d})
are given by the formula
\begin{eqnarray}
D_d(c|u) &=&
\displaystyle
\left(-\frac{c}{n}\right)^{n}
\sum_{|I|=d}
\sum_{J'\sqcup J''\subset I}
(-1)^{|J'|}
\frac{\theta^{(|J'|)}}{\theta}(u)
\frac{\partial^{J''}\Delta}{\Delta}(\lambda)
\cdot
\left(-\frac{n}{c}\partial\right)^{I\backslash(J'\sqcup J'')}
\nonumber
\\
&\stackrel{I\backslash J'=:I'}{=}&
\left(-\frac{c}{n}\right)^{n}
\sum_{k=0}^n
(-1)^k\frac{\theta^{(k)}}{\theta}(u)
\sum_{|I'|=d-k}
\sum_{J''\subset I'}
\frac{\partial^{J''}\Delta}{\Delta}(\lambda)
\cdot
\left(-\frac{n}{c}\partial\right)^{I'\backslash J''}.
\label{eq:coefofD(u)}
\end{eqnarray}
\hfill\fb
\label{prop:bibunop}
\end{prop}

Since the functions
$\{\frac{\theta^{(k)}}{\theta}(u)\}_{k=0,\cdots,n}$
are linearly independent,
we can consider the coefficient 
for $(-1)^k\frac{\theta^{(k)}}{\theta}(u)$
in $\left(-\frac{n}{c}\right)^n D_d(u)$
above.
 From $(\ref{eq:coefofD(u)})$,
it is $0$ for $d<k$, and
\begin{equation}
\sum_{|I|=d-k}
\sum_{J\subset I}
\frac{\partial^{J}\Delta}{\Delta}(\lambda)
\cdot
\left(-\frac{n}{c}\partial\right)^{I\backslash J}
=:D[d-k]
\label{eq:Debiard'sop}
\end{equation}
for $d\geq k$; it depends only on the difference $d-k$.
For example,
\begin{eqnarray*}
D[1]&=&
\sum_{i=1}^n(-\frac{n}{c}\partial_i)+\frac{\partial_i\Delta}{\Delta},\\
D[2]&=&
\sum_{i<j}(-\frac{n}{c}\partial_i)(-\frac{n}{c}\partial_j)
+\frac{\partial_i\Delta}{\Delta}(-\frac{n}{c}\partial_j)
+\frac{\partial_j\Delta}{\Delta}(-\frac{n}{c}\partial_i)
+\frac{\partial_i\partial_j\Delta}{\Delta}.
\end{eqnarray*}
These are the operators appeared in \cite{De} in the trigonometric case:

\medskip
{\bf Corollary to Proposition \ref{prop:bibunop}}
{\it
For $d=1,\cdots,n$, the operators
$D[d]$ (\ref{eq:Debiard'sop})
commute with each other.
}\hfill \fb
\medskip

To end this section, let us prove the previous formula (\ref{eq:limisOP}).
As is easily seen, the ``hamiltonian'' $H$
and these operators
are related by $H=D[1]^2-2D[2]$.
On the other hand,

\begin{prop}
$$
D[d-k]
=
{\left(-\frac{n}{c}\right)^n}
\sum_{k'=0}^d
(-1)^{d-k'}
{n-k' \choose n-d}
\frac{(-\frac{ck'}{n})^{k}}{k!}
\left.
\frac{\partial_\hbar^{d-k}\dot{M}_{k'}(c)}{(d-k)!}
\right|_{\hbar=0}.
$$
\end{prop}

{\bf Proof.}
By the defintion of $D[d-k]$, we have
$$
\sum_{d\geq k}
(-t)^{n-d}\frac{(-1)^k\theta^{(k)}}{\theta}(u)
D[d-k]
=
\left(-\frac{n}{c}\right)^{n}
\np
{\rm det}[\tilde{L}'_0-t]
\np.
$$
 From (\ref{eq:defD_d}) and Theorem \ref{thm:genfnc}, this is equal to
\begin{eqnarray*}
&&
{\left(-\frac{n}{c}\right)^n}
\lim_{\hbar\rightarrow 0}
\sum_{k'=0}^n
\frac{(-1-\hbar t)^{n-k'}}{\hbar^n}
\frac{\theta(u+\frac{ck'}{n}\hbar)}{\theta(u)}\dot{M}_{k'}(c)
\\
&=&
{\left(-\frac{n}{c}\right)^n}
\lim_{\hbar\rightarrow 0}
\sum_{k'=0}^n
(-1)^{n-k'}
\left(
\sum_{d=0}^{n-k'}
{n-k' \choose d}t^d\hbar^{d-n}
\right)
\left(
\sum_{k=0}^{\infty}
\frac{(\frac{ck'}{n}\hbar)^{k}}{k!}\frac{\theta^{(k)}}{\theta}(u)
\right)
\dot{M}_{k'}(c)
\\
&=&
{\left(-\frac{n}{c}\right)^n}
\sum_{k'=0}^n
(-1)^{n-k'}
\sum_{d=0}^{n-k'}
{n-k' \choose d}t^d
\sum_{k=0}^{n-d}
\frac{(\frac{ck'}{n})^{k}}{k!}\frac{\theta^{(k)}}{\theta}(u)
\left.
\frac{\partial_\hbar^{n-d-k}\dot{M}_{k'}(c)}{(n-d-k)!}
\right|_{\hbar=0}.
\end{eqnarray*}
Rearranging the last summation to obtain the
coefficient for
$(-t)^{n-d}\frac{(-1)^k\theta^{(k)}}{\theta}(u)$,
we get the assertion.
\hfill \fb

\medskip
Some particular cases of this proposition are:
\begin{eqnarray}
D[1]=D[1-0]&=&-\dot{M}_1
',\nonumber\\
D[2]=D[2-0]&=&
\frac{{n-2 \choose 0}\dot{M}_2
''-{n-1 \choose 1}\dot{M}_1
''}{2}
\label{eq:D[d]viaM_d}
=\frac{\dot{M}_2
''-(n-1)\dot{M}_1
''}{2},
\end{eqnarray}
where $'=\partial_\hbar|_{\hbar=0}$.
Now the formula (\ref{eq:limisOP}) can be verified as follows.
\begin{eqnarray}
\frac{\dot{M}_1^2-2\dot{M}_1-2\dot{M}_2+n}{\hbar^2}
&=&
(\frac{\dot{M}_1-n}{\hbar})^2
+2\frac{(n-1)(\dot{M}_1-n)-(\dot{M}_2-{n\choose 2})}{\hbar^2}
\nonumber
\\
&\stackrel{\hbar\rightarrow 0}{\longrightarrow}&
\dot{M}_1'^2+2\frac{(n-1)\dot{M}_1''-\dot{M}_2''}{2}
\stackrel{(\ref{eq:D[d]viaM_d})}{=}
D[1]^2-2D[2] = H.
\label{eq:l'Hospital}
\end{eqnarray}
In the last line we have used the l'Hospital rule
(since
$(n-1)\dot{M}_1'-\dot{M}_2'=0$,
we are to take the second derivative).

\label{sect:difflim}

\section{The invariant subspace spanned by theta functions}

In this section,
we would like to investigate a certain finite-dimensional
invariant subspace for $M_d(c|u)$ when $c$ is a nonnegative integer $l$.

Let us consider the space of level $l$ theta functions on the weight
space ${\bf h}^*$,
\begin{equation}
Th_l :=
\left\{
f:\mbox{\bf h}^* \stackrel{hol}{\rightarrow}\cpx
\vphantom{
\begin{array}{l}
f(\lambda + \alpha) = f(\lambda),\\
f(\lambda + \alpha\tau)
\end{array}
}
\right.
\left|
\begin{array}{l}
f(\lambda + \alpha) = f(\lambda),\\
f(\lambda + \alpha\tau)
 = f(\lambda)\exp [-2\pi il(\langle\lambda,\alpha\rangle
                           +\langle\alpha,\alpha\rangle\tau/2)]
\end{array}
(\forall \alpha \in Q)
\right\},
\label{eq:defTh_l}
\end{equation}
where $Q=Q(A_{n-1})$ stands for the root lattice,
and also let
$Th_l^{S(n)}$
be the subspace consisting of symmetric group action invariants.
This space is of dimension
$\frac{(l+n)!}{l!n!}$
and actually spanned by the level $l$ $A_{n-1}^{(1)}$ -characters.
For $j=0,\cdots,n-1$ let
$\Lambda_j$ be the classical part of the $j$-th fundamental
weight for the type $A_{n-1}^{(1)}$ root system
\cite{Kacbook}
and let
$$
\chi_{j}(\lambda)
:=
  \sum_{\nu \in \Lambda_j + Q}
   \exp 2\pi \sqrt{-1}
\left[\langle\lambda,\nu\rangle
+ \frac{\langle\nu,\nu\rangle\tau}{2}\right].
$$
Then more precisely we have
\begin{equation}
Th_l^{S(n)}
=
\bigoplus_{0\leq j_1\leq\cdots\leq j_l\leq n-1}
\cpx \chi_{j_1}\cdots\chi_{j_l}
{}.
\label{eq:basisofTh}
\end{equation}

In \cite{has93}, we have shown that $Th_l^{S(n)}$ is an
$A(R)$ submodule in $ {\cal O}({\bf h}^*)$:
\begin{thm}[\cite{has93}]
When the parameter $c$ is set to a nonnegative integer $l$,
we have
$$
L(l|u)^i_j Th_l^{S(n)} \subset Th_l^{S(n)}
$$
for any $u\in \cpx$ and $i,j=1,\cdots,n$ and therefore
$$
M_d(l|u) Th_l^{S(n)} \subset Th_l^{S(n)}.
$$
\end{thm}

(To be precise, the operators $L(l|u)^i_j$ as functions in $u$
have poles at $u\in {\bf Z} + \tau{\bf Z}$ and therefore the space
$Th_l^{S(n)}$ should be extended by tensoring ${\cal O}(\cpx)$, the
field of meromorphic functions in $u\in \cpx$.)
\hfill\fb

\medskip
Here we wish to establish a structure theorem for this
subspace as an $A(R)$- module.

To describe the result, let us give some definitions.
First we introduce an $A(R)$ module which is isomorphic to
the $A(R)$- comodule $V(\fb_{v_1}\otimes\cdots\otimes\fb_{v_l})$
as a vector space.
To avoid confusion let us denote  the space by
${\cal V}(\fb_{v_1}\otimes\cdots\otimes\fb_{v_l}).$
The $A(R)$ module structure for this space is given by \cite{FRT}
\begin{center}
$
L(u)^i_{i'}e^j
=
\sum_{j'} e^{j'} \cdot R(u-v)^{i,j}_{i',j'}
$
\quad
\mbox{\rm for $e^{j'} \in {\cal V}(\fb_v)$}
\end{center}
and in general we tensor this structure via the
comultiplication
\begin{equation}
L(u)^i_{i'}\mapsto
\sum_{i_1,\cdots,i_{l-1}}
L(u)^i_{i_1}
\otimes\cdots\otimes
L(u)^{i_{l-1}}_{i'}
\in
A(R)^{\otimes l}
\label{eq:coprod}
\end{equation}
so that we have
\begin{equation}
L(u)^i_{i'}{\rm v}
:=R(\fb_u,\fb_{v_1}\otimes\cdots\otimes\fb_{v_l})^i_{i'}
{\rm v}
\quad
\mbox{\rm for any}
\quad
{\rm v} \in {\cal V}(\fb_{v_1}\otimes\cdots\otimes\fb_{v_l}).
\label{eq:tensormod}
\end{equation}
Here 
$R(\fb_u,\fb_{v_1}\otimes\cdots\otimes\fb_{v_l})^i_{i'}$
stands for the operator on
$
{\cal V}(\fb_{v_1}\otimes\cdots\otimes\fb_{v_l})
$
 defined by
$R(\fb_u,\fb_{v_1}\otimes\cdots\otimes\fb_{v_l})
 e^i\otimes {\rm v}
=
\sum_{i'}
 e^{i'}\otimes
( R(\fb_u,\fb_{v_1}\otimes\cdots\otimes\fb_{v_l})^i_{i'}
 {\rm v}).$
Recall that
$
{\rm Im}(\check{R}(-\hbar))\subset \wedge^2(\cpx^n)
\subset
{\cal V}(\fb_v \otimes \fb_{v-\hbar})
$
(\ref{eq:imR=ext2}).
For more general lengthy case, we have
$$
{\cal I}
:=
\sum_{m=1}^{l-1}{\rm Im} \check{R}(-\hbar)^{k,k+1}
\subset
{\cal V}(\fb_v \otimes\cdots\otimes \fb_{v-(l-1)\hbar})
$$
and the subspace ${\cal I}$ is an $A(R)$-submodule.
This can be easily seen from the relation
\begin{eqnarray*}
&&
R(\fb_u,\fb_v
\otimes\cdots\otimes
\stackrel{(k-1)th}{\fb_{v-(k-1)\hbar}}
\otimes
\stackrel{k-th}{\fb_{v-k\hbar}}
\otimes\cdots\otimes
\fb_{v-(l-1)\hbar})^i_{i'}
 \check{R}(-\hbar)^{k,k+1}
\\
&=&
 \check{R}(-\hbar)^{k,k+1}
R(\fb_u,
\fb_v
\otimes\cdots\otimes
\stackrel{(k-1)th}{\fb_{v-k\hbar}}
\otimes
\stackrel{k-th}{\fb_{v-(k-1)\hbar}}
\otimes\cdots\otimes
 \fb_{-(l-1)\hbar})^i_{i'},
\end{eqnarray*}
which is of course a consequence of the YBE.
Now we define the quotient $A(R)$-module
$$
{\cal S}^l_v
:=
{\cal V}(\fb_v \otimes\cdots\otimes \fb_{v-(l-1)\hbar})
/\,{\cal I}.
$$
By the property (\ref{eq:imR=ext2}) we have
$$
\cdots\otimes(e^i\otimes e^j - e^j\otimes e^i)\otimes\cdots
\equiv 0 \,\,{\rm mod}\,\,{\cal I},
$$
and hence ${\cal S}^l_v$ is isomorphic to the usual
space of homogeneous polynomials of degree $l$
in $e^1, \cdots, e^n$:
\begin{equation}
{\cal S}^l_v
=
\bigoplus_{0\leq j_1\leq\cdots\leq j_l\leq n-1}
\cpx\, e^{j_1}\cdots e^{j_l},
\quad
e^{j_1}\cdots e^{j_l}
:=e^{j_1}\otimes\cdots\otimes e^{j_l} \,\,{\rm mod}\,\,{\cal I}.
\label{eq:basisofS}
\end{equation}

Comparing (\ref{eq:basisofTh}) and (\ref{eq:basisofS}),
they are obviously isomorphic as vector spaces by the
correspondence
\begin{equation}
\gamma:
e^{j_1}\cdots e^{j_l}
\mapsto
\chi_{j_1}\cdots\chi_{j_l}.
\label{eq:ThisomS}
\end{equation}

\begin{thm}
Let us modify the L-operator (\ref{eq:defofL}) in the case
$c=l$, a positive integer,
by the scalar function multiplication
\begin{equation}
L^\circ (l|u)^i_j:=
 \left(
  \prod_{s=0}^{l-1}
  \frac{\theta(u+s\hbar)}{\theta(\hbar)}
 \right)
\cdot L(l|u)^i_j
\label{eq:normalizedL}
\end{equation}
and consider the
$A(R)$-module structure on $Th^{S(n)}_l$
defined by $L(u)^i_j \mapsto L^\circ (l|u)^i_j$.
Then the correspondence (\ref{eq:ThisomS}) gives
the isomorphism of $A(R)$-modules
$$
\gamma :
{\cal S}^l_0
 \stackrel{\simeq}{\longrightarrow}
Th_l^{S(n)}
\quad
\mbox{
(the parameter $v$ for ${\cal S}^l_v$ is set to $0$)}.
$$
\label{thm:ThisomS}
\end{thm}

{\bf Remark.}
In the case $l=1$, we have shown in
\cite{has93}, Theorem 2(1)
that
\begin{equation}
L(l|u)^i_{i'}\chi_j
 =
 \frac{\theta(\hbar)}{\theta(u)}
  \sum_{j'=1}^n \chi_{j'} R(u)^{i\phantom{'}j}_{i'j'},
\label{eq:thminl=1}
\end{equation}
implying the isomorphism.
This immediately explains why the modification
(\ref{eq:normalizedL}) is natural.
 From (\ref{eq:thminl=1}), it follows directly that
$\chi_j (j=1,\cdots,n)$
are eigenfunctions for $M_1(1|u)$
(and actually 
for the system $\{M_d(1|u)\}_{d=1}^n$)
with the same eigenvalue
$\frac{\theta(\hbar)}{\theta(u)}
 \sum_{i=1}^n R(u)^{i0}_{i0}$.

\medskip
{\bf Proof of Theorem \ref{thm:ThisomS}.}
Based on the $l=1$ case (see the remark above),
the general case can be verified as follows.
\begin{eqnarray}
&&
\gamma
\left(
L(u)^i_{i'}(e^{j_1}\cdots e^{j_l})
\right)
(\lambda)
\nonumber
\\
&\stackrel{(\ref{eq:tensormod})}{=}&
\gamma
\left(R(\fb_u,\fb_0 \otimes\cdots\otimes \fb_{0-(l-1)\hbar})^i_{i'}
(e^{j_1}\otimes \cdots \otimes e^{j_l} \,{\rm mod}\,{\cal I})
\right)
(\lambda)
\nonumber
\\
&=&
\gamma
\left(
\sum_{i_1,\cdots,i_{l-1}=1}^n
(R(u+0)^i_{i_1} e^{j_1})
\otimes
(R(u+\hbar)^{i_1}_{i_2} e^{j_2})
\otimes\cdots\otimes
(R(u+(l-1)\hbar)^{i_{l-1}}_{i'} e^{j_l})
\,{\rm mod}\,{\cal I}
\right)
(\lambda)
\nonumber
\\
&=&
\gamma
\left(
\sum_{i_1,\cdots,i_{l-1}}
(R(u+0)^i_{i_1} e^{j_1})
(R(u+\hbar)^{i_1}_{i_2} e^{j_2})
\cdots
(R(u+(l-1)\hbar)^{i_{l-1}}_{i'} e^{j_l})
\right)
(\lambda)
\nonumber
\\
&=&
\left(
\sum_{i_1,\cdots,i_{l-1}}
(L^\circ(1|u+0)^i_{i_1} \chi_{j_1})
(L^\circ(1|u+\hbar)^{i_1}_{i_2} \chi_{j_2})
\cdots
(L^\circ(1|u+(l-1)\hbar)^{i_{l-1}}_{i'} \chi_{j_l})
\right)
(\lambda)
\label{eq:usel=1}
\\
&=&
\prod_{s=0}^{l-1}\frac{\theta(u+s\hbar)}{\theta(\hbar)}
\nonumber
\\
&&
\times
\sum_{i_1,\cdots,i_{l-1}=1}^n
\sum_{k_1,\cdots,k_l=1}^n
(
\iv     {u}
        {\lambda+\hbar\bar{\epsilon}_{k_1}}
        {\lambda}
        ^i
\ov     {u+\hbar}
        {\lambda+\hbar\bar{\epsilon}_{k_1}}
        {\lambda}
        _{i_1}
)
(
\iv     {u+\hbar}
        {\lambda+\hbar\bar{\epsilon}_{k_2}}
        {\lambda}
        ^{i_1}
\ov     {u+2\hbar}
        {\lambda+\hbar\bar{\epsilon}_{k_2}}
        {\lambda}
        _{i_2}
)
\nonumber
\\
&&
\phantom{
}
\cdots
(
\iv     {u+(l-1)\hbar}
        {\lambda+\hbar\bar{\epsilon}_{k_l}}
        {\lambda}
        ^{i_{l-1}}
\ov     {u+l\hbar}
        {\lambda+\hbar\bar{\epsilon}_{k_l}}
        {\lambda}
        _{i'}
)
\cdot
\chi_{j_1}(\lambda+\hbar\bar{\epsilon}_{k_1})
\chi_{j_2}(\lambda+\hbar\bar{\epsilon}_{k_2})
\cdots
\chi_{j_l}(\lambda+\hbar\bar{\epsilon}_{k_l})
\nonumber
\\
&\stackrel{(\ref{eq:dualityofitv})}{=}&
 \prod_{s=0}^{l-1}\frac{\theta(u+s\hbar)}{\theta(\hbar)}
\nonumber\\&&
\times
\sum_{k_1,\cdots,k_l=1}^n
\,
\iv     {u}
        {\lambda+\hbar\bar{\epsilon}_{k_1}}
        {\lambda}
        ^i
\delta_{k_1,k_2}\cdots\delta_{k_{l-1},k_l}
\ov     {u+l\hbar}
        {\lambda+\hbar\bar{\epsilon}_{k_l}}
        {\lambda}
        _{i'}
\,\cdot\,
\chi_{j_1}(\lambda+\hbar\bar{\epsilon}_{k_1})
\cdots
\chi_{j_l}(\lambda+\hbar\bar{\epsilon}_{k_l})
\label{eq:key!}
\\
&=&
 \prod_{s=0}^{l-1}\frac{\theta(u+s\hbar)}{\theta(\hbar)}
\cdot
\sum_{k=1}^n
\iv     {u}
        {\lambda+\hbar\bar{\epsilon}_{k}}
        {\lambda}
        ^i
\ov     {u+l\hbar}
        {\lambda+\hbar\bar{\epsilon}_{k}}
        {\lambda}
        _{i'}
\left(
\chi_{j_1}\chi_{j_2}\cdots\chi_{j_l}
\right)
(\lambda+\hbar\bar{\epsilon}_{k})
\nonumber
\\
&=&
\left(
L^\circ(l|u)^i_{i'}
(\chi_{j_1}\cdots\chi_{j_l})
\right)(\lambda)
=
\left(
L^\circ(l|u)^i_{i'}
\gamma(e^{j_1}\cdots e^{j_l})
\right)(\lambda).
\nonumber
\end{eqnarray}

\begin{oekaki}\refstepcounter{oekaki}\addtocounter{oekaki}{-1}
\label{fig:freponTh}\caption{}\end{oekaki}

\noindent
We used the $l=1$ case in (\ref{eq:usel=1}).
The line (\ref{eq:key!}) is the most essential step,
where as shown in Fig.\ref{fig:freponTh}
the situation allows us to
employ the duality relation for the intertwining vectors
(\ref{eq:dualityofitv}).
\hfill \fb

\medskip
Thus there is an interesting ``representation theoretic" invariant
subspaces for our operators.
This space would be identified with the space of Weyl group invariant
theta functions in \cite{EK2}, where they considered an affine
analogue of Sutherland operator and its diagonalization.
Note that,
since our operators recover the elliptic Calogero-Moser system
(\ref{eq:limisOP}),
our system can be regarded as the q-analogue of the critical level case of
Etingof-Kirilov's theory.
Thus the existence of
$n$ 
commuting operators
$M_k(c|u)$ is in a coincidence with the
large center of the (completed) quantum affine enveloping algebra
at the critical level.


\section{Discussion}
In this paper we realized Ruijsenaars' commuting difference system as
commuting transfer matrices in the solvable lattice model point of
view.
In other word, this system can be considered as a function space
realization (a Schr\"{o}dinger picture) of the
one-dimensional spin chain defined by the elliptic R-matrix as its
Boltzmann weight.
Regarded as a lattice model, this system is very trivial because it
contains only one site of freedom.
It is quite natural to think of the multi-sites case as well,
that is,
consider the (fused) L-operators and their traces
acting on the tensor space
$
Th^{S(n)}_{l_1} \otimes \cdots \otimes Th^{S(n)}_{l_N}
$.
On this function space we can consider the action of the commuting
transfer matrices again.
The resulting system can be regarded as a difference version
of the elliptic curve case of generalized Gaudin models of type A
:\cite{ER}, \cite{Frenkel}, \cite{FFR}, \cite{GN}, \cite{Ne}.
In these papers Gaudin Hamiltonians are generalized
upon a geometric setting,
and the commuting Hamiltonians can be obtained by
the pole expansion of
the characteristic polynomial of the (quantized) Higgs field.
Therefore the factorized L-operator
or the difference Lax matrix $\tilde{L}$ (Proposition 3)
can be said as a difference analogue of the Higgs field.
This suggests that the intertwining vectors,
the building block for the factorized L-operator,
would also allow some intrinsic definition
(geometric as well as representation theoretic).
The paper \cite{Kr-Zab} is in this direction with
the geometric treatment of the intertwining vectors and
consequently the factorized L-operator for Baxter's R-matrix.
But still it seems that the origin of the weight space, or the
dynamical parameter $\lambda$, needs to be elucidated.
Another related work is \cite{C95}, where Cherednik announced the
construction of elliptic difference operators based on his Hecke
algebra representation technique. Whether and how we can find the
relationship with his thoeory and the present construction would be
also an important problem.


\medskip
{\bf Acknowledgement.}
The author expresses his sincere gratitude to
Prof. E. K. Sklyanin, Prof. L. D. Faddeev,
Prof. E. Date, Prof. R. Hotta and Prof. A. Tsuchiya
for their kind interest and encouragement.
Also thanks are due to
Prof. M. Noumi, Prof. H. Awata, Prof. G. Kuroki, Prof. T. Takebe,
Prof. H. Ochiai and Prof. van Diejen
for fruitful discussions and comments.
During the preparation of this manuscript the author visited
City University, London
and
Universit\'{e} Louis Pasteur, Strasbourg,
respectively supported by the JSPS and by the Kawai foundation.
It is his great pleasure to thank for the hospitality at these universities,
especially for Prof. P. P. Martin and Prof. G. Schiffmann.

\appendix

\section{Proof of Lemma 1.}
Recall $Y_{r<s}:=1$ iff $r<s$ and $0$ otherwise, and put
$$
A^{(d)}_{s,s'} =
A^{(d)}(u,\lambda,\mu)_{s,s'}:=
 \prod_{r=1}^d
  \theta
    \left(\mu_r - \lambda_{s'} +\hbar Y_{r<s} + \delta_{r,s}(u-(s-1)\hbar)
    \right).
$$
Then
$A^{(d)}= \left[A^{(d)}_{s,s'}\right]_{s,s'=1,\cdots,d}$
is the matrix in problem.
Lemma \ref{lem:qFay} (\ref{eq:qFay})
trivially holds for the case $d=1$ ; suppose that $d-1$ case is true.

\medskip
First we let $\mu_d=\lambda_{s'}$ for some $s'=1,\cdots,d$,
and observe
\begin{itemize}
\item
$
A^{(d)}_{s,s'}=0 (s\leq d), \;
A^{(d)}_{d,s'}
=\theta(u-(d-1)\hbar)\prod_{r=1}^{d-1}\theta(\mu_r-\lambda_{s'}+\hbar)
$
\, and
\item
$
A^{(d)}_{s,t}=A^{(d-1)}_{s,t}\cdot\theta(\lambda_{s'}-\lambda_s)
\;\; (1\leq s,t \leq d-1).
$
\end{itemize}
Then the determinant is simply a product of
$A^{(d)}_{d,s'}$
and the remaining $(d-1)\times(d-1)$ determinant:
$$
{\rm det}\left[A^{(d)}_{s,t}\right]_{s,t=1,\cdots,d}
=
{\rm det}\left[A^{(d-1)}_{s,t}\right]_{s=1,\cdots,d-1;t\neq s'}
\cdot
\theta(u-(d-1)\hbar)
\prod_{r=1}^{d-1}\theta(\mu_r-\lambda_{s'}+\hbar)
\prod_{s,\neq s'}\theta(\lambda_{s'}-\lambda_s)
$$
By the induction hypothesis, the right hand side is equal to
$$
\theta(u+\sum_{r=1}^d(\mu_r-\lambda_r))
\prod_{s=1}^{d-1}\theta(u-s\hbar)
\prod_{1\leq s<s'\leq d}
 \theta(\lambda_{s'}-\lambda_{s})
 \theta(\hbar + \mu_s-\mu_{s'})
$$
as desired.
That is, the formula in question holds for $d$ special values
$\mu_d=\lambda_{s'} (s'=1,\cdots,d)$ of $\mu_d$.

\medskip
A standard complex analysis tells that,
for an entire nonzero holomprphic function $f(z)$ with the properties
$$
f(z+1)=\exp[-2\pi\sqrt{-1}B]f(z), \,
f(z+\tau)=\exp[-2\pi\sqrt{-1}(Cz+C')]f(z)
$$
have $C$ zeros $z_1,\cdots,z_{C}$
with multiplicities inside a fundamental period parallelogram
and
$$
\sum_{j=1}^{C}{z_j}
\equiv
B\tau+\frac{C}{2}-C'\;
{\rm mod}\; {\bf Z}+{\bf  Z}\tau$$
holds.

\medskip
Consider the difference
$$
{\rm det}\left[A^{(d)}_{s,t}\right]_{s,t=1,\cdots,d}
-\theta(u+\sum_{r=1}^d(\mu_r-\lambda_r))
\prod_{s=1}^{d-1}\theta(u-s\hbar)
\prod_{1\leq s<s'\leq d}
 \theta(\lambda_{s'}-\lambda_{s})
 \theta(\hbar + \mu_s-\mu_{s'})
$$
as a function in $\mu_d=z$ and denote it by $f(z)$.
Then it is easily checked that
$f(z+1)=(-1)^d f(z)$ and
\begin{eqnarray*}
f(z+\tau)
&=&
\exp[2\pi\sqrt{-1}(-u+(d-1)\hbar)]\prod_{s=1}^d
(-\exp[2\pi\sqrt{-1}(-\mu_d+\lambda_s-\tau/2)])\cdot f(z)
\\
&=&
\exp\left[
-2\pi\sqrt{-1}
 (dz-\sum_{s=1}^d\lambda_s+d\frac{\tau+1}{2}+u-(d-1)\hbar)
\right]
\cdot f(z)
{}.
\end{eqnarray*}
Therefore, if $f$ is not identically zero,
we should have $d$ zeros modulo the period lattice whose summation
should be given by
$$
\frac{d}{2}\tau+\frac{d}{2}-(-\sum_{s=1}^d\lambda_s+d\frac{\tau+1}{2}+u-(d-1)\hbar)
\equiv\sum_{s=1}^d\lambda_s-u+(d-1)\hbar
\; {\rm mod} \,{\bf Z}+{\bf Z}\tau.
$$
On the other hand, the previous observation shows that we have $d$
zeros $\lambda_1,\cdots,\lambda_d$ whose summation is
$
\sum_{s=1}^d\lambda_s \not\equiv \sum_{s=1}^d\lambda_s-u+(d-1)\hbar
$.
Since $u$ and $\hbar$ are arbitrary,
this is a contradiction.
Thus $f$ should be identically zero, proving the lemma.
\hfill\fb





\newpage
\begin{center}{\bf\huge Figure Captions}\end{center}
\medskip

\unitlength=0.6pt
\begin{center}
\begin{picture}(292.00,167.00)(100.00,406.00)
\put(274.00,556.00){\vector(1,-1){130.00}}
\put(400.00,559.00){\vector(-1,-1){130.00}}
\put(297.00,563.00){\vector(1,-1){23.00}}
\put(378.00,563.00){\vector(-1,-1){20.00}}
\put(150.00,496.00)
{\makebox(0,0)[cc]
{\large $R(u-v)^{i,j}_{i',j'}=\check{R}(u-v)^{i,j}_{j',i'}=$}}
\put(266.00,573.00){\makebox(0,0)[cc]{\large $i$}}
\put(403.00,573.00){\makebox(0,0)[cc]{\large $j$}}
\put(261.00,415.00){\makebox(0,0)[cc]{\large $j'$}}
\put(407.00,412.00){\makebox(0,0)[cc]{\large $i'$}}
\put(334.00,527.00){\makebox(0,0)[cc]{\large $u-v$}}
\put(310.00,566.00){\makebox(0,0)[cc]{\large $u$}}
\put(358.00,566.00){\makebox(0,0)[cc]{\large $v$}}
\end{picture}

\medskip
{\large Fig.\ref{oe:fRijij}. 
Matrix elements of the R-matrix.}
\end{center}

\vfil
\vskip 5mm

\begin{minipage}{7.0cm}
\begin{center}{

\unitlength=0.7pt
\begin{picture}(280.00,102.00)(121.00,418.00)
\put(401.00,519.00){\vector(-1,0){162.00}}
\put(320.00,518.00){\vector(0,-1){76.00}}
\put(200.00,470.00){\makebox(0,0)[rc]
{\large $\phi(u)^\mu_{\lambda,j} =$}}
\put(225.00,520.00){\makebox(0,0)[cc]{\large $\lambda$}}
\put(425.00,520.00){\makebox(0,0)[cc]{\large $\mu$}}
\put(334.00,498.00){\makebox(0,0)[cc]{\large $u$}}
\put(321.00,428.00){\makebox(0,0)[cc]{\large $j$}}
\end{picture}

\vskip 1ex

{\large Fig.{\ref{oe:foutitv}}.
\\
The outgoing intertwining vector.}

}\end{center}
\end{minipage}
\phantom{XXXXX}
%
%
\begin{minipage}{7.5cm}
\begin{center}{
\unitlength=0.7pt
\begin{picture}(289.00,96.00)(123.00,440.00)
\put(321.00,521.00){\vector(0,-1){78.00}}
\put(401.00,441.00){\vector(-1,0){157.00}}
\put(175.00,480.00){\makebox(0,0)[cc]
{\large $\bar{\phi}(u)^{\mu,i}_\lambda = $}}
\put(225.00,440.00){\makebox(0,0)[cc]{\large $\lambda$}}
\put(412.00,440.00){\makebox(0,0)[cc]{\large $\mu$}}
\put(339.00,461.00){\makebox(0,0)[cc]{\large $u$}}
\put(321.00,536.00){\makebox(0,0)[cc]{\large $j$}}
\end{picture}

\vskip 2ex

{\large Fig.{\ref{oe:fincitv}}.
\\
The incoming intertwining vector.}

}\end{center}
\end{minipage}
\vskip 5mm
\vfil
\begin{center}
\unitlength=0.9pt
\begin{picture}(462.00,142.00)(60.00,368.00)
\put(400.00,440.00){\line(1,0){104.00}}
\put(400.00,430.00){\line(1,0){99.00}}
\put(225.00,490.00){\line(-1,0){74.00}}
\put(225.00,390.00){\line(-1,0){74.00}}
\put(125.00,465.00){\line(0,-1){49.00}}
\put(175.00,490.00){\vector(0,-1){46.00}}
\put(175.00,444.00){\vector(0,-1){55.00}}
\put(450.00,490.00){\vector(0,-1){49.00}}
\put(450.00,430.00){\vector(0,-1){49.00}}
\put(161.00,440.00){\makebox(0,0)[cc]{\large $k$}}
\put(177.00,510.00){\makebox(0,0)[cc]{\large $u$}}
\put(176.00,379.00){\makebox(0,0)[cc]{\large $u$}}
\put(467.00,455.00){\makebox(0,0)[cc]{\large $u$}}
\put(466.00,415.00){\makebox(0,0)[cc]{\large $u$}}
\put(512.00,433.00){\makebox(0,0)[cc]{\large $\lambda$}}
\put(392.00,431.00){\makebox(0,0)[cc]{\large $\mu$}}
\put(90.00,440.00){\makebox(0,0)[cc]{\large $\displaystyle\sum_k$}}
\put(261.00,439.00){\makebox(0,0)[cc]{\large $=\delta_{j,j'}$}}
\put(552.00,439.00){\makebox(0,0)[cc]{\large $=\delta_{k,k'}$}}
\put(452.00,509.00){\makebox(0,0)[cc]{\large $k$}}
\put(451.00,368.00){\makebox(0,0)[cc]{\large $k'$}}
\put(358.00,434.00){\makebox(0,0)[cc]{\large $\displaystyle\sum_\lambda$}}
\put(227.00,501.00){\makebox(0,0)[cc]{\large $\lambda+\hbar{\bar\epsilon}_j$}}
\put(228.00,377.00){\makebox(0,0)[cc]{\large
$\lambda+\hbar{\bar\epsilon}_{j'}$}}
\put(135.00,446.00){\makebox(0,0)[cc]{\large $\lambda$}}
\bezier50(150.00,490.00)(125.00,490.00)(125.00,465.00)
\bezier50(125.00,415.00)(125.00,390.00)(150.00,390.00)
\end{picture}
\\
\vskip 5mm
{\large
Fig.\ref{fig:dualityofitv}
\\
The duality relation.
}
\end{center}

\vfil
\begin{center}
\unitlength=0.8pt
\begin{picture}(359.00,178.00)(100.00,393.00)
\put(321.00,563.00){\vector(0,-1){75.00}}
\put(402.00,486.00){\vector(-1,0){160.00}}
\put(402.00,476.00){\vector(-1,0){160.00}}
\put(321.00,476.00){\vector(0,-1){70.00}}
\put(140.00,479.00){\makebox(0,0)[cc]
{
{\large $L(c|u)^i_j f(\mu) =$}
{\large $\displaystyle{\sum_{k=1}^n}$}
$\left( \rule[-10ex]{0pt}{20ex}\right.$
}}
\put(232.00,480.00){\makebox(0,0)[cc]{\large $\mu$}}
\put(321.00,571.00){\makebox(0,0)[cc]{\large $i$}}
\put(321.00,393.00){\makebox(0,0)[cc]{\large $j$}}
\put(339.00,502.00){\makebox(0,0)[cc]{\large $u$}}
\put(335.00,461.00){\makebox(0,0)[lc]{\large $u+c\hbar$}}
\put(415.00,483.00){\makebox(0,0)[lc]
{{\large $\mu+\hbar\bar{\epsilon}_k$}
$\left. \rule[-10ex]{0pt}{20ex}\right)$
}}
\put(500.00,480.00){\makebox(0,0)[lc]{\large $f(\mu+\hbar\bar{\epsilon}_k)$}}
\end{picture}

\medskip
{\large Fig.{\ref{fig:ffactorL}}
\\
The factorized L-operator.}
\end{center}
\vfil

\begin{center}
\unitlength=0.7pt
\begin{picture}(290.00,191.00)(120.00,383.00)
\put(360.00,560.00){\vector(-1,-1){80.00}}
\put(441.00,480.00){\vector(-1,1){80.00}}
\put(441.00,479.00){\vector(-1,-1){81.00}}
\put(360.00,398.00){\vector(-1,1){79.00}}
\put(159.00,479.00){\makebox(0,0)[cc]
{\large
$
\check{W}
\left[\begin{array}{ccc}
                        & \mu&\\
                \lambda &u&\nu\\
                        & \mu'& \\\end{array}\right]
=$}}
\put(263.00,481.00){\makebox(0,0)[cc]{\large $\lambda$}}
\put(362.00,574.00){\makebox(0,0)[cc]{\large $\mu$}}
\put(361.00,383.00){\makebox(0,0)[cc]{\large $\mu'$}}
\put(449.00,480.00){\makebox(0,0)[cc]{\large $\nu$}}
\put(361.00,479.00){\makebox(0,0)[cc]{\large $u$}}
\end{picture}

{\large Fig.{\ref{fig:ffacewt}}
\\
The face Boltzmann weight.}
\end{center}
\vfil

\begin{center}
\unitlength=1.0pt
\begin{picture}(364.00,168.00)(117.00,411.00)
\put(197.00,561.00){\vector(-1,0){76.00}}
\put(279.00,560.00){\vector(-1,0){73.00}}
\put(161.00,559.00){\vector(0,-1){32.00}}
\put(242.00,558.00){\vector(0,-1){31.00}}
\put(243.00,521.00){\vector(-1,-1){82.00}}
\put(163.00,519.00){\vector(1,-1){80.00}}
\put(401.00,561.00){\vector(-1,-1){61.00}}
\put(400.00,440.00){\vector(-1,1){60.00}}
\put(461.00,502.00){\vector(-1,1){59.00}}
\put(462.00,502.00){\vector(-1,-1){59.00}}
\put(396.00,431.00){\vector(-1,1){59.00}}
\put(467.00,494.00){\vector(-1,-1){59.00}}
\put(369.00,458.00){\vector(-1,-1){33.00}}
\put(438.00,463.00){\vector(1,-1){32.00}}
\put(401.00,434.00){\circle*{8}}
\put(162.00,522.00){\circle*{6}}
\put(243.00,521.00){\circle*{6}}
\put(117.00,573.00){\makebox(0,0)[cc]{\large $\lambda$}}
\put(200.00,576.00){\makebox(0,0)[cc]{\large $\mu$}}
\put(279.00,576.00){\makebox(0,0)[cc]{\large $\nu$}}
\put(161.00,575.00){\makebox(0,0)[cc]{\large $u$}}
\put(243.00,576.00){\makebox(0,0)[cc]{\large $v$}}
\put(152.00,432.00){\makebox(0,0)[cc]{\large $j$}}
\put(243.00,427.00){\makebox(0,0)[cc]{\large $i$}}
\put(199.00,506.00){\makebox(0,0)[cc]{\large $u-v$}}
\put(299.00,499.00){\makebox(0,0)[cc]{\large $=$}}
\put(324.00,501.00){\makebox(0,0)[cc]{\large $\lambda$}}
\put(401.00,500.00){\makebox(0,0)[cc]{\large $u-v$}}
\put(399.00,579.00){\makebox(0,0)[cc]{\large $\mu$}}
\put(481.00,500.00){\makebox(0,0)[cc]{\large $\nu$}}
\put(332.00,415.00){\makebox(0,0)[cc]{\large $j$}}
\put(471.00,415.00){\makebox(0,0)[cc]{\large $i$}}
\put(368.00,430.00){\makebox(0,0)[cc]{\large $v$}}
\put(437.00,430.00){\makebox(0,0)[cc]{\large $u$}}
\end{picture}

{\large Fig.{\ref{fig:fovprop}}.
\\
The intertwining property of the outgoing vectors.
\\
Dots represent summation indices.}
\end{center}
\vfil
\begin{center}
\unitlength=1.0pt
\begin{picture}(379.00,171.00)(116.00,318.00)
\put(158.00,483.00){\vector(1,-1){76.00}}
\put(239.00,401.00){\vector(0,-1){36.00}}
\put(238.00,482.00){\vector(-1,-1){79.00}}
\put(159.00,403.00){\vector(0,-1){38.00}}
\put(279.00,362.00){\vector(-1,0){74.00}}
\put(193.00,362.00){\vector(-1,0){73.00}}
\put(419.00,462.00){\vector(-1,-1){60.00}}
\put(416.00,342.00){\vector(-1,1){57.00}}
\put(477.00,403.00){\vector(-1,1){56.00}}
\put(479.00,402.00){\vector(-1,-1){61.00}}
\put(416.00,467.00){\vector(-1,-1){61.00}}
\put(484.00,405.00){\vector(-1,1){62.00}}
\put(354.00,471.00){\vector(1,-1){28.00}}
\put(483.00,466.00){\vector(-1,-1){28.00}}
\put(159.00,403.00){\circle*{6}}
\put(238.00,403.00){\circle*{6}}
\put(418.00,465.00){\circle*{8}}
\put(148.00,497.00){\makebox(0,0)[cc]{\large $i$}}
\put(240.00,500.00){\makebox(0,0)[cc]{\large $j$}}
\put(198.00,468.00){\makebox(0,0)[cc]{\large $u-v$}}
\put(114.00,355.00){\makebox(0,0)[cc]{\large $\lambda$}}
\put(200.00,355.00){\makebox(0,0)[cc]{\large $\mu$}}
\put(286.00,356.00){\makebox(0,0)[cc]{\large $\nu$}}
\put(159.00,348.00){\makebox(0,0)[cc]{\large $v$}}
\put(240.00,349.00){\makebox(0,0)[cc]{\large $u$}}
\put(300.00,422.00){\makebox(0,0)[cc]{\large $=$}}
\put(418.00,403.00){\makebox(0,0)[cc]{\large $u-v$}}
\put(347.00,489.00){\makebox(0,0)[cc]{\large $i$}}
\put(485.00,486.00){\makebox(0,0)[cc]{\large $j$}}
\put(381.00,467.00){\makebox(0,0)[cc]{\large $u$}}
\put(458.00,468.00){\makebox(0,0)[cc]{\large $v$}}
\put(343.00,405.00){\makebox(0,0)[cc]{\large $\lambda$}}
\put(419.00,329.00){\makebox(0,0)[cc]{\large $\mu$}}
\put(493.00,401.00){\makebox(0,0)[cc]{\large $\nu$}}

\end{picture}

{\large Fig.{\ref{fig:fivprop}}.
\\
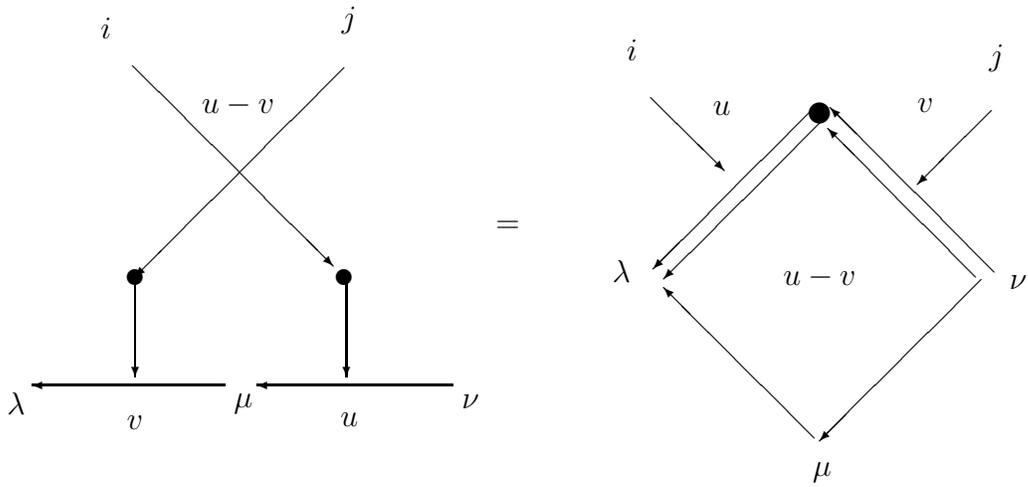
The intertwining property of the incoming vectors.
}
\end{center}

\vfil
\begin{center}
\unitlength=0.9pt
\begin{picture}(508.00,153.00)(62.00,404.00)
\put(90.00,479.00){\vector(1,-1){62.00}}
\put(151.00,480.00){\vector(-1,-1){62.00}}
\put(90.00,540.00){\vector(0,-1){29.00}}
\put(119.00,512.00){\vector(-1,0){60.00}}
\put(120.00,507.00){\vector(-1,0){61.00}}
\put(90.00,506.00){\vector(0,-1){25.00}}
\put(150.00,541.00){\vector(0,-1){30.00}}
\put(180.00,511.00){\vector(-1,0){59.00}}
\put(180.00,507.00){\vector(-1,0){59.00}}
\put(151.00,506.00){\vector(0,-1){24.00}}
\put(300.00,524.00){\vector(-1,-1){45.00}}
\put(345.00,480.00){\vector(-1,1){42.00}}
\put(346.00,479.00){\vector(-1,-1){46.00}}
\put(300.00,433.00){\vector(-1,1){45.00}}
\put(348.00,484.00){\vector(-1,1){44.00}}
\put(352.00,529.00){\vector(-1,-1){25.00}}
\put(295.00,528.00){\vector(-1,-1){43.00}}
\put(250.00,531.00){\vector(1,-1){25.00}}
\put(349.00,474.00){\vector(-1,-1){47.00}}
\put(328.00,452.00){\vector(1,-1){23.00}}
\put(295.00,429.00){\vector(-1,1){44.00}}
\put(272.00,451.00){\vector(-1,-1){26.00}}
\put(451.00,540.00){\vector(1,-1){60.00}}
\put(510.00,540.00){\vector(-1,-1){62.00}}
\put(511.00,478.00){\vector(0,-1){27.00}}
\put(540.00,451.00){\vector(-1,0){61.00}}
\put(540.00,447.00){\vector(-1,0){61.00}}
\put(511.00,446.00){\vector(0,-1){25.00}}
\put(478.00,451.00){\vector(-1,0){59.00}}
\put(480.00,446.00){\vector(-1,0){60.00}}
\put(449.00,479.00){\vector(0,-1){26.00}}
\put(449.00,446.00){\vector(0,-1){26.00}}
\put(122.00,509.00){\circle*{6}}
\put(90.00,481.00){\circle*{6}}
\put(152.00,479.00){\circle*{6}}
\put(300.00,528.00){\circle*{6}}
\put(300.00,428.00){\circle*{6}}
\put(449.00,478.00){\circle*{6}}
\put(512.00,478.00){\circle*{6}}
\put(479.00,448.00){\circle*{6}}
\put(210.00,480.00){\makebox(0,0)[cc]{\large $=$}}
\put(390.00,484.00){\makebox(0,0)[cc]{\large $=$}}
\put(90.00,556.00){\makebox(0,0)[cc]{\large $i$}}
\put(151.00,557.00){\makebox(0,0)[cc]{\large $j$}}
\put(88.00,404.00){\makebox(0,0)[cc]{\large $j"$}}
\put(153.00,405.00){\makebox(0,0)[cc]{\large $i"$}}
\put(47.00,510.00){\makebox(0,0)[cc]{\large $\lambda$}}
\put(244.00,481.00){\makebox(0,0)[cc]{\large $\lambda$}}
\put(408.00,447.00){\makebox(0,0)[cc]{\large $\lambda$}}
\put(195.00,509.00){\makebox(0,0)[cc]{\large $\nu$}}
\put(364.00,480.00){\makebox(0,0)[cc]{\large $\nu$}}
\put(550.00,447.00){\makebox(0,0)[cc]{\large $\nu$}}
\put(242.00,541.00){\makebox(0,0)[cc]{\large $i$}}
\put(449.00,553.00){\makebox(0,0)[cc]{\large $i$}}
\put(363.00,543.00){\makebox(0,0)[cc]{\large $j$}}
\put(511.00,556.00){\makebox(0,0)[cc]{\large $j$}}
\put(243.00,404.00){\makebox(0,0)[cc]{\large $j"$}}
\put(450.00,404.00){\makebox(0,0)[cc]{\large $j"$}}
\put(361.00,408.00){\makebox(0,0)[cc]{\large $i"$}}
\put(512.00,408.00){\makebox(0,0)[cc]{\large $i"$}}
\put(102.00,524.00){\makebox(0,0)[cc]{\large $u$}}
\put(164.00,527.00){\makebox(0,0)[cc]{\large $v$}}
\put(347.00,504.00){\makebox(0,0)[cc]{\large $v$}}
\put(462.00,463.00){\makebox(0,0)[cc]{\large $v$}}
\put(274.00,525.00){\makebox(0,0)[cc]{\large $u$}}
\put(527.00,464.00){\makebox(0,0)[cc]{\large $u$}}
\put(120.00,468.00){\makebox(0,0)[cc]{\large $u-v$}}
\put(302.00,479.00){\makebox(0,0)[cc]{\large $u-v$}}
\put(480.00,529.00){\makebox(0,0)[cc]{\large $u-v$}}
\put(94.00,497.00){\makebox(0,0)[lc]{\large $u+c\hbar$}}
\put(350.00,451.00){\makebox(0,0)[cc]{\large $u+c\hbar$}}
\put(517.00,430.00){\makebox(0,0)[lc]{\large $u+c\hbar$}}
\put(157.00,490.00){\makebox(0,0)[lc]{\large $v+c\hbar$}}
\put(271.00,428.00){\makebox(0,0)[cc]{\large $v+c\hbar$}}
\put(454.00,430.00){\makebox(0,0)[lc]{\large $v+c\hbar$}}
\end{picture}
\vskip 1cm
{\large
Fig.\ref{fig:fRLL_LLR}
\\
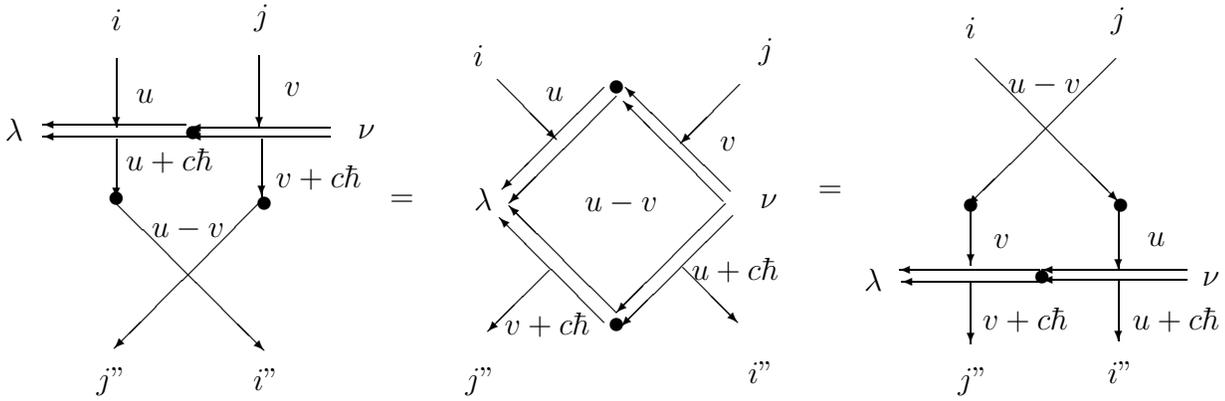
Proof of the RLL=LLR relation.
}
\end{center}
\vfil

\begin{center}
\unitlength=1.0pt
\begin{picture}(189.00,183.00)(150.00,438.00)
\put(441.00,599.00){\vector(-1,-1){159.00}}
\put(401.00,599.00){\vector(1,-1){33.00}}
\put(363.00,559.00){\vector(1,-1){38.00}}
\put(281.00,480.00){\vector(1,-1){42.00}}
\put(346.00,522.00){\circle*{2}}
\put(325.00,501.00){\circle*{2}}
\put(337.00,512.00){\circle*{2}}
\put(279.00,502.00){\makebox(0,0)[cc]{\large $\fb_{u_1}$}}
\put(359.00,577.00){\makebox(0,0)[cc]{\large $\fb_{u_{k-1}}$}}
\put(400.00,615.00){\makebox(0,0)[cc]{\large $\fb_{u_{k}}$}}
\put(448.00,616.00){\makebox(0,0)[cc]{\large $\fb_{v}$}}
\put(259.00,539.00)
{\makebox(0,0)[rc]{\large
$\check{R}(\fb_{u_1}\otimes\cdots\otimes\fb_{u_k},\fb_v) =$
}}
\end{picture}
\vskip 1cm
{\large
Fig.\ref{fig:ffusedR1}.
\\
The matrix $\check{R}(\fb_{u_1}\otimes\cdots\otimes\fb_{u_k},\fb_v)$.
}
\end{center}
\vfil

\begin{center}
\unitlength=1.0pt
\begin{picture}(175.00,166.00)(153.00,419.00)
\put(272.00,510.00){\vector(1,-1){90.00}}
\put(291,530){\vector(1,-1){90.00}}
\put(330.00,569.00){\vector(1,-1){93.00}}
\put(361.00,569.00){\vector(-1,-1){91.00}}
\put(381,550){\vector(-1,-1){91.00}}
\put(422.00,510.00){\vector(-1,-1){91.00}}
\put(311.00,532.00){\circle*{2}}
\put(316.5,537.5){\circle*{2}}
\put(322.00,543.00){\circle*{2}}

\put(383.00,531.00){\circle*{2}}
\put(388.5,526.5){\circle*{2}}
\put(394.00,521.00){\circle*{2}}

\put(368.00,493.00){\makebox(0,0)[cc]{\LARGE $\cdots$}}

\put(265.00,523.00){\makebox(0,0)[cc]{\large $\fb_{u_1}$}}
\put(285.00,543.00){\makebox(0,0)[cc]{\large $\fb_{u_2}$}}
\put(320.00,581.00){\makebox(0,0)[cc]{\large $\fb_{u_k}$}}
\put(368.00,585.00){\makebox(0,0)[cc]{\large $\fb_{v_1}$}}
\put(388.00,565.00){\makebox(0,0)[cc]{\large $\fb_{v_2}$}}
\put(428.00,525.00){\makebox(0,0)[cc]{\large $\fb_{v_l}$}}
\put(253.00,493.00){\makebox(0,0)[rc]{\large
$\check{R}
(\fb_{u_1}\otimes\cdots\otimes\fb_{u_k},
\fb_{v_1}\otimes\cdots\otimes\fb_{v_l}) =$}}
\end{picture}
\vskip 1cm
{\large
Fig.\ref{fig:ffusedR2}
\\
The fused R-matrix.
}\end{center}
\vfil

\begin{center}
\unitlength=1.0pt
\begin{picture}(152.00,150.00)(260.00,440.00)
\put(400.00,558.00){\vector(-1,-1){118.00}}
\put(281.00,558.00){\vector(1,-1){118.00}}
\put(400.00,521.00){\line(0,-1){6.00}}
\put(400.00,515.00){\vector(-1,-1){76.00}}
\put(337.00,542.00){\makebox(0,0)[cc]{\large $\cdots$}}
\put(361.00,440.00){\makebox(0,0)[cc]{\large $\cdots$}}
\put(267.00,572.00){\makebox(0,0)[cc]{\large $\fb_{u-(k-1)\hbar}$}}
\put(352.00,573.00){\makebox(0,0)[cc]{\large $\fb_{u-\hbar}$}}
\put(404.00,573.00){\makebox(0,0)[cc]{\large $\fb_u$}}
\put(256.00,500.00){\makebox(0,0)[rc]{\large $\pi_{1^k}=$}}
\put(362.00,559.00){\line(1,-1){38}}
\end{picture}

\vskip 1cm
{\large
Fig.\ref{fig:fvproj}.
\\
The fusion operator $\pi_{1^k}$.
}
\end{center}
\vskip 5mm
\vfil

\begin{center}
\unitlength=1.0pt
\begin{picture}(250.00,14.00)(221.00,479.00)
\put(401.00,479.00){\vector(-1,0){80.00}}
\put(250.00,481.00){\makebox(0,0)[cc]{\large
$e_\lambda^{\lambda+\hbar\bar{\epsilon}_i}$ =}}
\put(308.00,480.00){\makebox(0,0)[cc]{\large $\lambda$}}
\put(424.00,479.00){\makebox(0,0)[cc]{\large $\lambda+\hbar\bar{\epsilon}_i$}}
\put(363.00,493.00){\makebox(0,0)[cc]{\large $u$}}
\end{picture}

\vskip 1cm
{\large
Fig.\ref{fig:f1steppath}.
\\
One step path.
}\end{center}
\vskip 5mm
\vfil

\begin{center}
\unitlength=1.0pt
\begin{picture}(370.00,33.00)(119.00,464.00)
\put(402.00,480.00){\vector(-1,0){60.00}}
\put(340.00,480.00){\vector(-1,0){60.00}}
\put(500.00,480.00){\vector(-1,0){60.00}}
\put(422.00,480.00){\makebox(0,0)[cc]{\large $\cdots$}}
\put(250.00,480.00){\makebox(0,0)[rc]
{\large
${
P}_\lambda^\nu
=\sum_{\mu_1,\cdots,\mu_{k-1}}{\bf C}$}}
\put(280.00,465.00){\makebox(0,0)[cc]{\large $\lambda$}}
\put(345.00,465.00){\makebox(0,0)[cc]{\large $\mu_1$}}
\put(401.00,465.00){\makebox(0,0)[cc]{\large $\mu_{2}$}}
\put(445.00,465.00){\makebox(0,0)[cc]{\large $\mu_{k-1}$}}
\put(500.00,465.00){\makebox(0,0)[cc]{\large $\nu$}}
\put(319.00,490.00){\makebox(0,0)[cc]{\large $u_1$}}
\put(374.00,490.00){\makebox(0,0)[cc]{\large $u_2$}}
\put(477.00,490.00){\makebox(0,0)[cc]{\large $u_k$}}
\end{picture}
\vskip 1cm

  {\large
Fig.\ref{fig:flongpath}.
\\
   Space of paths.}
\end{center}

\vskip 5mm
\vfil
\unitlength=1.0pt
\begin{picture}(353.00,214.00)(100.00,402.00)
\put(271.00,509.00){\vector(-1,1){60.00}}
\put(211.00,569.00){\vector(-1,-1){60.00}}
\put(271.00,509.00){\vector(-1,-1){61.00}}
\put(210.00,448.00){\vector(-1,1){59.00}}
\put(271.00,599.00){\vector(-1,0){61.00}}
\put(210.00,599.00){\vector(-1,0){60.00}}
\put(481.00,509.00){\vector(-1,1){60.00}}
\put(421.00,569.00){\vector(-1,-1){60.00}}
\put(480.00,509.00){\vector(-1,-1){61.00}}
\put(419.00,448.00){\vector(-1,1){59.00}}
\put(481.00,420.00){\vector(-1,0){59.00}}
\put(422.00,420.00){\vector(-1,0){62.00}}
\put(213.00,582.00){\circle*{6}}
\put(300.00,510.00){\makebox(0,0)[cc]{\large $:=\sum_{\mu'}$}}
\put(211.00,510.00){\makebox(0,0)[cc]{\large $\check{W}(\fb_u,\fb_v)$}}
\put(142.00,600.00){\makebox(0,0)[cc]{\large $\lambda$}}
\put(211.00,614.00){\makebox(0,0)[cc]{\large $\mu$}}
\put(271.00,616.00){\makebox(0,0)[cc]{\large $\nu$}}
\put(347.00,511.00){\makebox(0,0)[cc]{\large $\lambda$}}
\put(420.00,584.00){\makebox(0,0)[cc]{\large $\mu$}}
\put(495.00,512.00){\makebox(0,0)[cc]{\large $\nu$}}
\put(421.00,433.00){\makebox(0,0)[cc]{\large $\mu'$}}
\put(488.00,435.00){\makebox(0,0)[cc]{\large $\nu'$}}
\put(352.00,434.00){\makebox(0,0)[cc]{\large $\lambda'$}}
\put(421.00,508.00){\makebox(0,0)[cc]{\large $u-v$}}
\end{picture}

\medskip
\begin{center}
{\large
Fig.\ref{fig:ffaceop}
\\
Definition of the face operator $\check{W}(\fb_u,\fb_v)$.}
\end{center}
\vfil

\begin{center}
\unitlength=1.0pt
\begin{picture}(319.00,163.00)(115.00,399.00)
\put(199.00,562.00){\vector(-1,-1){40.00}}
\put(239.00,521.00){\vector(-1,1){40.00}}
\put(239.00,521.00){\vector(-1,-1){39.00}}
\put(200.00,482.00){\vector(-1,1){41.00}}
\put(159.00,523.00){\vector(-1,-1){39.00}}
\put(239.00,442.00){\vector(-1,1){40.00}}
\put(199.00,482.00){\vector(-1,-1){39.00}}
\put(160.00,443.00){\vector(-1,1){39.00}}
\put(239.00,441.00){\vector(-1,-1){42.00}}
\put(199.00,401.00){\vector(-1,1){38.00}}
\put(439.00,481.00){\vector(-1,1){40.00}}
\put(399.00,521.00){\vector(-1,1){41.00}}
\put(358.00,562.00){\vector(-1,-1){38.00}}
\put(399.00,521.00){\vector(-1,-1){40.00}}
\put(359.00,481.00){\vector(-1,1){39.00}}
\put(439.00,481.00){\vector(-1,-1){41.00}}
\put(398.00,440.00){\vector(-1,-1){39.00}}
\put(359.00,401.00){\vector(-1,1){41.00}}
\put(398.00,442.00){\vector(-1,1){41.00}}
\put(357.00,483.00){\vector(-1,-1){39.00}}
\put(275.00,479.00){\makebox(0,0)[cc]{\large $=$}}
\put(201.00,520.00){\makebox(0,0)[cc]{\large $\check{W}(\fb_v,\fb_w)$}}
\put(160.00,479.00){\makebox(0,0)[cc]{\large $\check{W}(\fb_u,\fb_w)$}}
\put(201.00,440.00){\makebox(0,0)[cc]{\large $\check{W}(\fb_u,\fb_v)$}}
\put(359.00,520.00){\makebox(0,0)[cc]{\large $\check{W}(\fb_u,\fb_v)$}}
\put(361.00,441.00){\makebox(0,0)[cc]{\large $\check{W}(\fb_v,\fb_w)$}}
\put(401.00,480.00){\makebox(0,0)[cc]{\large $\check{W}(\fb_u,\fb_w)$}}
\end{picture}

{\large
Fig.\ref{fig:ffaceYBE}.
\\
Face Yang-Baxter equation.
}\end{center}

\vfil

\begin{center}
\unitlength=1.0pt
\begin{picture}(242.00,320.00)(150.00,319.00)
\put(400.00,639.00){\vector(-1,-1){39.00}}
\put(441.00,599.00){\vector(-1,1){39.00}}
\put(441.00,599.00){\vector(-1,-1){40.00}}
\put(401.00,559.00){\vector(-1,1){40.00}}
\put(440.00,521.00){\vector(-1,1){38.00}}
\put(402.00,559.00){\vector(-1,-1){41.00}}
\put(441.00,518.00){\vector(-1,-1){39.00}}
\put(402.00,479.00){\vector(-1,1){40.00}}
\put(440.00,360.00){\vector(-1,1){39.00}}
\put(401.00,399.00){\vector(-1,-1){40.00}}
\put(442.00,359.00){\vector(-1,-1){40.00}}
\put(402.00,319.00){\vector(-1,1){40.00}}
\put(359.00,599.00){\vector(-1,-1){119.00}}
\put(360.00,359.00){\vector(-1,1){118.00}}
\put(401.00,400.00){\vector(-1,1){120.00}}
\put(362.00,517.00){\vector(-1,-1){77.00}}
\put(362.00,520.00){\vector(-1,1){39.00}}
\put(400.00,440.00){\circle*{3}}
\put(360.00,480.00){\circle*{3}}
\put(380.00,460.00){\circle*{3}}
\put(320.00,440.00){\circle*{3}}
\put(310.00,450.00){\circle*{3}}
\put(330.00,430.00){\circle*{3}}
\put(200.00,480.00){\makebox(0,0)[rc]{\large $\Pi_{1^k}:=$}}
\put(402.00,599.00){\makebox(0,0)[cc]{\large $-\hbar$}}
\put(401.00,519.00){\makebox(0,0)[cc]{\large $-\hbar$}}
\put(402.00,359.00){\makebox(0,0)[cc]{\large $-\hbar$}}
\put(361.00,558.00){\makebox(0,0)[cc]{\large $-2\hbar$}}
\put(282.00,478.00){\makebox(0,0)[cc]{\large $-(k-1)\hbar$}}
\put(361.00,397.00){\makebox(0,0)[cc]{\large $-2\hbar$}}
\end{picture}

{\large
Fig.\ref{fig:ffproj}.
\\
Face fusion operator.}
\end{center}
\vfil

\begin{center}
\unitlength=1.0pt
\begin{picture}(280.00,124.00)(161.00,491.00)
\put(242.00,560.00){\vector(-1,0){81.00}}
\put(201.00,558.00){\vector(0,-1){35.00}}
\put(242.00,599.00){\vector(-1,0){81.00}}
\put(441.00,558.00){\vector(-1,0){81.00}}
\put(400.00,556.00){\vector(0,-1){35.00}}
\put(200.00,578.00){\circle*{4}}
\put(161.00,614.00){\makebox(0,0)[cc]{\large $\lambda$}}
\put(243.00,615.00){\makebox(0,0)[cc]{\large $\mu$}}
\put(216.00,538.00){\makebox(0,0)[cc]{\large $\fb_u$}}
\put(280.00,559.00){\makebox(0,0)[cc]{\large $:=$}}
\put(320.00,560.00){\makebox(0,0)[cc]{\large $\sum_j$}}
\put(361.00,573.00){\makebox(0,0)[cc]{\large $\lambda$}}
\put(440.00,574.00){\makebox(0,0)[cc]{\large $\mu$}}
\put(413.00,542.00){\makebox(0,0)[cc]{\large $u$}}
\put(401.00,510.00){\makebox(0,0)[cc]{\large $j$}}
\put(401.00,485.00){\makebox(0,0)[cc]{\large $e^j$}}
\end{picture}

{\large
Fig.\ref{fig:fov}.
\\
Definition of $\phi(\fb_u)$.
}\end{center}
\vfil

\begin{center}
\unitlength=1.0pt
\begin{picture}(285.00,82.00)(159.00,460.00)
\put(200.00,518.00){\vector(0,-1){38.00}}
\put(241.00,479.00){\vector(-1,0){82.00}}
\put(401.00,519.00){\vector(0,-1){39.00}}
\put(440.00,479.00){\vector(-1,0){80.00}}
\put(440.00,460.00){\vector(-1,0){79.00}}
\put(201.00,538.00){\makebox(0,0)[cc]{\large $e^j$}}
\put(200.00,528){\circle*{5}}
\put(215.00,496.00){\makebox(0,0)[cc]{\large $\fb_u$}}
\put(280.00,479.00){\makebox(0,0)[cc]{\large $:=$}}
\put(313.00,483.00){\makebox(0,0)[cc]{\large $\sum_{\lambda,k}$}}
\put(352.00,468.00){\makebox(0,0)[cc]{\large $\lambda$}}
\put(444.00,468.00){\makebox(0,0)[lc]{\large $\lambda+\hbar\bar{\epsilon}_k$}}
\put(416.00,501.00){\makebox(0,0)[cc]{\large $u$}}
\put(402.00,542.00){\makebox(0,0)[cc]{\large $j$}}
\end{picture}

\vskip 1cm

{\large
Fig.\ref{fig:fiv}.
\\
Definition of $\bar{\phi}(\fb_u)$.
}\end{center}
\vfil

\begin{center}
\unitlength=1.0pt
\begin{picture}(391.00,284.00)(89.00,396.00)
\put(381.00,546.00){\line(1,-1){66.00}}
\put(210.00,659.00){\vector(-1,-1){28.00}}
\put(182.00,631.00){\vector(-1,-1){30.00}}
\put(152.00,601.00){\vector(-1,-1){61.00}}
\put(242.00,629.00){\vector(-1,1){32.00}}
\put(242.00,629.00){\vector(-1,-1){120.00}}
\put(241.00,570.00){\vector(-1,1){60.00}}
\put(241.00,569.00){\vector(-1,-1){90.00}}
\put(210.00,419.00){\vector(-1,1){121.00}}
\put(240.00,450.00){\vector(-1,1){120.00}}
\put(241.00,449.00){\vector(-1,-1){30.00}}
\put(241.00,419.00){\vector(-1,0){30.00}}
\put(211.00,419.00){\vector(-1,0){32.00}}
\put(179.00,419.00){\vector(-1,0){30.00}}
\put(149.00,419.00){\vector(-1,0){28.00}}
\put(121.00,419.00){\vector(-1,0){31.00}}
\put(109.00,418.00){\vector(0,-1){22.00}}
\put(141.00,417.00){\vector(0,-1){20.00}}
\put(228.00,418.00){\vector(0,-1){19.00}}
\put(480.00,569.00){\vector(-1,0){30.00}}
\put(450.00,569.00){\vector(-1,0){32.00}}
\put(418.00,569.00){\vector(-1,0){30.00}}
\put(388.00,569.00){\vector(-1,0){27.00}}
\put(361.00,569.00){\vector(-1,0){30.00}}
\put(350.00,567.00){\vector(0,-1){18.00}}
\put(379.00,568.00){\vector(0,-1){20.00}}
\put(468.00,567.00){\vector(0,-1){21.00}}
\put(468.00,546.00){\vector(-1,-1){123.00}}
\put(351.00,547.00){\vector(1,-1){122.00}}
\put(442.00,442.00){\vector(-1,-1){18.00}}
\put(211.00,510.00){\circle*{2}}
\put(211.00,520.00){\circle*{2}}
\put(211.00,500.00){\circle*{2}}
\put(183.00,402.00){\makebox(0,0)[cc]{\large $\cdots$}}
\put(181.00,539.00){\makebox(0,0)[cc]{\large $\Pi_{1^k}$}}
\put(286.00,537.00){\makebox(0,0)[cc]{\large $=$}}
\put(391.00,429.00){\makebox(0,0)[cc]{\large $\cdots$}}
\put(422.00,545.00){\makebox(0,0)[cc]{\large $\cdots$}}
\put(374.00,486.00){\makebox(0,0)[cc]{\large $\pi_{1^k}$}}
\put(109.00,380.00){\makebox(0,0)[cc]{\large $u$}}
\put(142.00,380.00){\makebox(0,0)[cc]{\large $u-\hbar$}}
\put(227.00,380.00){\makebox(0,0)[cc]{\large $u-(k-1)\hbar$}}
\put(211.00,629.00){\makebox(0,0)[cc]{\large $-\hbar$}}
\put(211.00,568.00){\makebox(0,0)[cc]{\large $-\hbar$}}
\put(211.00,449.00){\makebox(0,0)[cc]{\large $-\hbar$}}
\put(181.00,599.00){\makebox(0,0)[cc]{\large $-2\hbar$}}
\put(121.00,540.00){\makebox(0,0)[cc]{\large $-(k-1)\hbar$}}
\put(150.00,508.00){\makebox(0,0)[cc]{\large $-(k-2)\hbar$}}
\put(349.00,590.00){\makebox(0,0)[cc]{\large $u-(k-1)\hbar$}}
\put(425.00,589.00){\makebox(0,0)[cc]{\large $\cdots$}}
\put(467.00,590.00){\makebox(0,0)[cc]{\large $u$}}
\bezier62(449.00,479.00)(466.00,459.00)(437.00,437.00)
\end{picture}
\end{center}
\vskip 5mm
\begin{center}
{\large
Fig.\ref{fig:piphi=phiPi}
\\
The relation (\ref{eq:piphi=phiPi}).
}
\end{center}
\vfil

\begin{center}
\unitlength=1.0pt
\begin{picture}(436.00,309.00)(60.00,340.00)
\put(201.00,640.00){\vector(0,-1){37.00}}
\put(240.00,600.00){\vector(-1,0){79.00}}
\put(240.00,576.00){\vector(-1,0){79.00}}
\put(200.00,574.00){\vector(0,-1){34.00}}
\put(401.00,637.00){\vector(0,-1){40.00}}
\put(400.00,570.00){\vector(0,-1){35.00}}
\put(361.00,599.00){\vector(0,-1){11.00}}
\put(361.00,571.00){\vector(0,1){15.00}}
\put(163.00,439.00){\vector(0,-1){38.00}}
\put(163.00,401.00){\vector(0,-1){39.00}}
\put(361.00,438.00){\vector(0,-1){34.00}}
\put(361.00,404.00){\vector(0,-1){41.00}}
\put(242.00,439.00){\vector(0,-1){38.00}}
\put(242.00,401.00){\vector(0,-1){37.00}}
\put(121.00,417.00){\vector(0,-1){17.00}}
\put(121.00,380.00){\vector(0,1){16.00}}
\put(201.00,587.00){\circle{6}}
\put(163.00,400.00){\circle*{5}}
\put(242.00,400.00){\circle*{5}}
\put(361.00,399.00){\circle*{5}}
\put(400.50,586.50){\oval(79.00,103.00)}
\put(274.50,399.50){\oval(307.00,79.00)}
\put(202.00,649.00){\makebox(0,0)[cc]{\large $e^I$}}
\put(199.00,527.00){\makebox(0,0)[cc]{\large $e_I$}}
\put(141.00,588.00){\makebox(0,0)[cc]{\large $\lambda$}}
\put(252.00,586.00){\makebox(0,0)[cc]{\large $\lambda+\hbar\bar{\epsilon}_J$}}
\put(105.00,592.00){\makebox(0,0)[cc]{\large $\sum_I$}}
\put(284.00,583.00){\makebox(0,0)[cc]{\large $=$}}
\put(399.00,581.00){\makebox(0,0)[cc]{\large $V^{\otimes k}$}}
\put(342.00,585.00){\makebox(0,0)[cc]{\large $\lambda$}}
\put(480.00,584.00){\makebox(0,0)[cc]{\large $\lambda+\hbar\bar{\epsilon}_J$}}
\put(216.00,611.00){\makebox(0,0)[cc]{\large $1^k_u$}}
\put(413.00,553.00){\makebox(0,0)[cc]{\large $1^k_u$}}
\put(220.00,559.00){\makebox(0,0)[cc]{\large $1^k_{u+c\hbar}$}}
\put(420.00,620.00){\makebox(0,0)[cc]{\large $1^k_{u+c\hbar}$}}
\put(62.00,400.00){\makebox(0,0)[cc]
{\large $=\sum_{\sigma\in S(d)}{\rm sgn}(\sigma)$}}
\put(108.00,400.00){\makebox(0,0)[cc]{\large $\lambda$}}
\put(304.00,400.00){\makebox(0,0)[cc]{\large $\cdots$}}
\put(475.00,402.00){\makebox(0,0)[cc]{\large $\lambda+\hbar\bar{\epsilon}_J$}}
\put(175.00,379.00){\makebox(0,0)[lc]{\large $\fb_u$}}
\put(254.00,379.00){\makebox(0,0)[lc]{\large $\fb_{u-\hbar}$}}
\put(367.00,379.00){\makebox(0,0)[lc]{\large $\fb_{u-(d-1)\hbar}$}}
\put(175.00,422.00){\makebox(0,0)[lc]{\large $\fb_{u+c\hbar}$}}
\put(254.00,422.00){\makebox(0,0)[lc]{\large $\fb_{u-\hbar+c\hbar}$}}
\put(367.00,422.00){\makebox(0,0)[lc]{\large $\fb_{u-(d-1)\hbar+c\hbar}$}}
\put(164.00,340.00){\makebox(0,0)[cc]{\large $+\hbar\bar{\epsilon}_{j_1}$}}
\put(243.00,340.00){\makebox(0,0)[cc]{\large $+\hbar\bar{\epsilon}_{j_2}$}}
\put(362.00,340.00){\makebox(0,0)[cc]{\large $+\hbar\bar{\epsilon}_{j_d}$}}
\put(165.00,458.00){\makebox(0,0)[cc]{\large
$+\hbar\bar{\epsilon}_{j_\sigma(1)}$}}
\put(242.00,459.00){\makebox(0,0)[cc]{\large
$+\hbar\bar{\epsilon}_{j_\sigma(2)}$}}
\put(362.00,460.00){\makebox(0,0)[cc]{\large
$+\hbar\bar{\epsilon}_{j_\sigma(d)}$}}
\end{picture}

\medskip

{\large
Fig.\ref{fig:fcalctr}
\\
Calculation of the trace.
}\end{center}

\vfil

\begin{center}
\unitlength=1.0pt
\begin{picture}(240.00,82.00)(200.00,479.00)
\put(362.00,559.00){\vector(-1,0){82.00}}
\put(321.00,559.00){\vector(0,-1){38.00}}
\put(321.00,521.00){\vector(0,-1){40.00}}
\put(362.00,480.00){\vector(-1,0){81.00}}
\put(321.00,519.00){\circle*{6}}
\put(200.00,520.00){\makebox(0,0)[cc]{\large $\tilde{L}(c|u)^i_j =$}}
\put(266.00,559.00){\makebox(0,0)[cc]{\large $\lambda$}}
\put(266.00,479.00){\makebox(0,0)[cc]{\large $\lambda$}}
\put(385.00,559.00){\makebox(0,0)[cc]{\large $\lambda+\hbar\bar{\epsilon}_i$}}
\put(385.00,479.00){\makebox(0,0)[cc]{\large $\lambda+\hbar\bar{\epsilon}_j$}}
\put(325,495){\makebox(0,0)[lc]{\large $u$}}
\put(325,545){\makebox(0,0)[lc]{\large $u+c\hbar$}}

\put(410.00,522.00){\makebox(0,0)[cc]{\large $T^\hbar_i$}}
\end{picture}

\vskip 1cm
{\large
Fig.\ref{fig:fLtilde}.
\\
The matrix $\tilde{L}(c|u)$.
}
\end{center}
\vfil

\begin{center}
\unitlength=1.0pt
\begin{picture}(310.00,306.00)(143.00,347.00)
\put(200.00,639.00){\vector(0,-1){38.00}}
\put(240.00,601.00){\vector(-1,0){82.00}}
\put(240.00,596.00){\vector(-1,0){82.00}}
\put(200.00,596.00){\vector(0,-1){34.00}}
\put(200.00,560.00){\vector(0,-1){36.00}}
\put(240.00,522.00){\vector(-1,0){82.00}}
\put(240.00,517.00){\vector(-1,0){82.00}}
\put(200.00,516.00){\vector(0,-1){35.00}}
\put(200.00,444.00){\vector(0,-1){38.00}}
\put(240.00,405.00){\vector(-1,0){82.00}}
\put(240.00,400.00){\vector(-1,0){82.00}}
\put(200.00,398.00){\vector(0,-1){37.00}}
\put(401.00,558.00){\vector(0,-1){37.00}}
\put(441.00,521.00){\vector(-1,0){81.00}}
\put(441.00,516.00){\vector(-1,0){81.00}}
\put(401.00,515.00){\vector(0,-1){34.00}}
\put(200.00,469.00){\circle*{4}}
\put(200.00,461.00){\circle*{4}}
\put(200.00,453.00){\circle*{4}}
\put(200.00,559.00){\circle*{6}}
\put(200.00,480.00){\circle*{6}}
\put(200.00,443.00){\circle*{6}}
\put(206.00,612.00){\makebox(0,0)[lc]{\large $u$}}
\put(206.00,581.00){\makebox(0,0)[lc]{\large $u+\hbar$}}
\put(206.00,537.00){\makebox(0,0)[lc]{\large $u+\hbar$}}
\put(206.00,501.00){\makebox(0,0)[lc]{\large $u+2\hbar$}}
\put(206.00,422.00){\makebox(0,0)[lc]{\large $u+(l-1)\hbar$}}
\put(206.00,386.00){\makebox(0,0)[lc]{\large $u+l\hbar$}}
\put(143.00,600.00){\makebox(0,0)[cc]{\large $\lambda$}}
\put(143.00,518.00){\makebox(0,0)[cc]{\large $\lambda$}}
\put(146.00,401.00){\makebox(0,0)[cc]{\large $\lambda$}}
\put(253.00,603.00){\makebox(0,0)[lc]{\large
$\lambda+\hbar\bar{\epsilon}_{k_1}$}}
\put(252.00,519.00){\makebox(0,0)[lc]{\large
$\lambda+\hbar\bar{\epsilon}_{k_2}$}}
\put(252.00,404.00){\makebox(0,0)[lc]{\large
$\lambda+\hbar\bar{\epsilon}_{k_3}$}}
\put(200.00,653.00){\makebox(0,0)[cc]{\large $i$}}
\put(199.00,347.00){\makebox(0,0)[cc]{\large $j$}}
\put(313.00,498.00){\makebox(0,0)[cc]{\large $=$}}
\put(346.00,520.00){\makebox(0,0)[cc]{\large $\lambda$}}
\put(468.00,520.00){\makebox(0,0)[cc]{\large $\lambda+\hbar\bar{\epsilon}_k$}}
\put(401.00,571.00){\makebox(0,0)[cc]{\large $i$}}
\put(403.00,468.00){\makebox(0,0)[cc]{\large $j$}}
\put(410.00,535.00){\makebox(0,0)[lc]{\large $u$}}
\put(410.00,500.00){\makebox(0,0)[lc]{\large $u+l\hbar$}}
\put(360.00,426.00){\makebox(0,0)[lc]
{\large $
\times\delta^k_{k_1}\delta^{k_1}_{k_2}\cdots\delta^{k_{l-1}}_{k_l}$
}}
\end{picture}

{\large
Fig.\ref{fig:freponTh}.
\\
Proof of Theorem 4.
}\end{center}


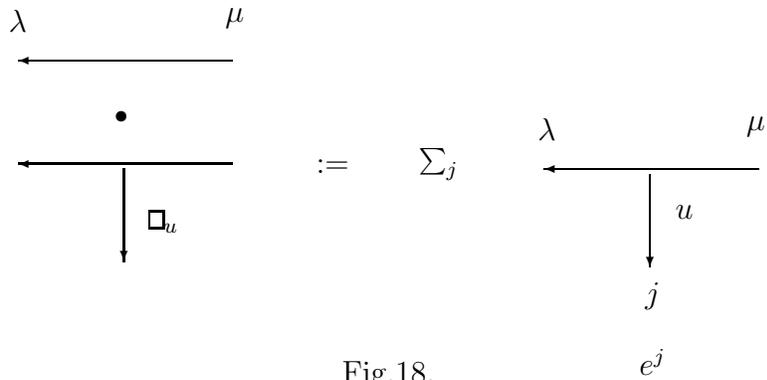
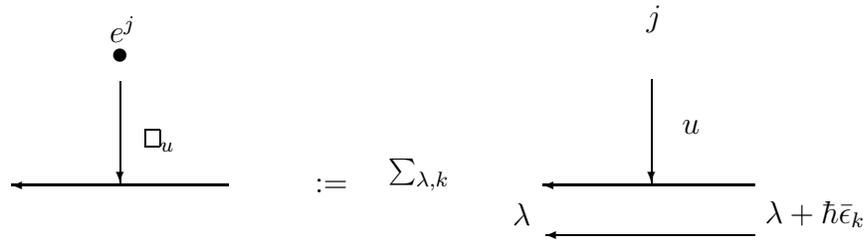
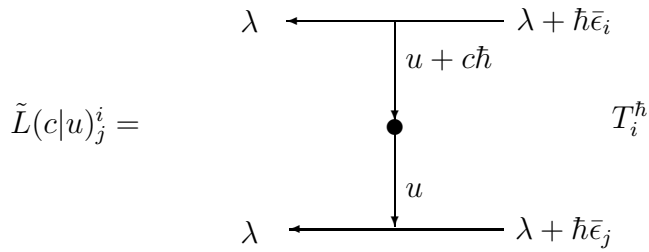
\end{document}